\documentclass[11pt]{article}
\usepackage[centertags]{amsmath}
\usepackage{amsfonts}
\usepackage{amssymb}
\usepackage{amsthm}
\usepackage{newlfont}
\usepackage{fullpage}%,backref}
\usepackage[colorlinks,linkcolor=blue, pdfstartview=FitH]{hyperref}
\theoremstyle{plain}
\newtheorem{thm}{Theorem}[section]
\newtheorem{theorem}[thm]{Theorem}
\newtheorem{lemma}[thm]{Lemma}
\newtheorem{corollary}[thm]{Corollary}
\newtheorem{proposition}[thm]{Proposition}

\theoremstyle{definition}

\newtheorem{remark}[thm]{Remark}

\newtheorem{definition}[thm]{Definition}

\newtheorem{example}[thm]{Example}

\newtheorem{defn-thm}[thm]{Definition-Theorem}

\numberwithin{equation}{section}

\begin{document}
%\begin{CJK*}{GBK}{song}

% ------------------------------------------------------------------------
%%%%%%%%%%%%%%%%%%%%%%%%%%%%%% New command%%%%%%%%%%%%%%%%%%%%%%%%%%%%%%%%%%%%%%%%%%%%

\newcommand{\sA}{{\mathcal A}}
\newcommand{\sB}{{\mathcal B}}
\newcommand{\sC}{{\mathcal C}}
\newcommand{\sD}{{\mathcal D}}
\newcommand{\sE}{{\mathcal E}}
\newcommand{\sF}{{\mathcal F}}
\newcommand{\sG}{{\mathcal G}}
\newcommand{\sH}{{\mathcal H}}
\newcommand{\sI}{{\mathcal I}}
\newcommand{\sJ}{{\mathcal J}}
\newcommand{\sK}{{\mathcal K}}
\newcommand{\sL}{{\mathcal L}}
\newcommand{\sM}{{\mathcal M}}
\newcommand{\sN}{{\mathcal N}}
\newcommand{\sO}{{\mathcal O}}
\newcommand{\sP}{{\mathcal P}}
\newcommand{\sQ}{{\mathcal Q}}
\newcommand{\sR}{{\mathcal R}}
\newcommand{\sS}{{\mathcal S}}
\newcommand{\sT}{{\mathcal T}}
\newcommand{\sU}{{\mathcal U}}
\newcommand{\sV}{{\mathcal V}}
\newcommand{\sW}{{\mathcal W}}
\newcommand{\sX}{{\mathcal X}}
\newcommand{\sY}{{\mathcal Y}}
\newcommand{\sZ}{{\mathcal Z}}
%mathfrak
\newcommand{\ssA}{{\mathfrak A}}
\newcommand{\ssB}{{\mathfrak B}}
\newcommand{\ssC}{{\mathfrak C}}
\newcommand{\ssD}{{\mathfrak D}}
\newcommand{\ssE}{{\mathfrak E}}
\newcommand{\ssF}{{\mathfrak F}}
\newcommand{\ssg}{{\mathfrak g}}
\newcommand{\ssH}{{\mathfrak H}}
\newcommand{\ssI}{{\mathfrak I}}
\newcommand{\ssJ}{{\mathfrak J}}
\newcommand{\ssK}{{\mathfrak K}}
\newcommand{\ssL}{{\mathfrak L}}
\newcommand{\ssM}{{\mathfrak M}}
\newcommand{\ssN}{{\mathfrak N}}
\newcommand{\ssO}{{\mathfrak O}}
\newcommand{\ssP}{{\mathfrak P}}
\newcommand{\ssQ}{{\mathfrak Q}}
\newcommand{\ssR}{{\mathfrak R}}
\newcommand{\ssS}{{\mathfrak S}}
\newcommand{\ssT}{{\mathfrak T}}
\newcommand{\ssU}{{\mathfrak U}}
\newcommand{\ssV}{{\mathfrak V}}
\newcommand{\ssW}{{\mathfrak W}}
\newcommand{\ssX}{{\mathfrak X}}
\newcommand{\ssY}{{\mathfrak Y}}
\newcommand{\ssZ}{{\mathfrak Z}}

\newcommand{\B}{{\mathbb B}}
\newcommand{\C}{{\mathbb C}}
\newcommand{\D}{{\mathbb D}}
\newcommand{\E}{{\mathbb E}}
\newcommand{\F}{{\mathbb F}}
\newcommand{\G}{{\mathbb G}}
\renewcommand{\H}{{\mathbb H}}
\newcommand{\J}{{\mathbb J}}
\newcommand{\K}{{\mathbb K}}
\renewcommand{\L}{{\mathbb L}}
\newcommand{\M}{{\mathbb M}}
\newcommand{\N}{{\mathbb N}}
\renewcommand{\P}{{\mathbb P}}
\newcommand{\Q}{{\mathbb Q}}
\newcommand{\R}{{\mathbb R}}

\newcommand{\T}{{\mathbb T}}
\newcommand{\U}{{\mathbb U}}
\newcommand{\V}{{\mathbb V}}
\newcommand{\W}{{\mathbb W}}
\renewcommand{\S}{{\mathbb S}}
\newcommand{\Y}{{\mathbb Y}}
\newcommand{\Z}{{\mathbb Z}}
\newcommand{\id}{{  id}}

\newcommand{\A}{\"{a}}
\newcommand{\rank}{{rank}}
\newcommand{\END}{{\mathbb E}{ nd}}
\newcommand{\End}{{ End}}
\newcommand{\Hg}{{Hg}}
\newcommand{\tr}{{ tr}}
\newcommand{\Tr}{{ Tr}}
\newcommand{\SL}{{ SL}}
\newcommand{\PSL}{{ PSL}}
\newcommand{\Gl}{{ Gl}}
\newcommand{\Cor}{{ Cor}}
\newcommand{\Gal}{{ Gal}}
\newcommand{\GL}{{ GL}}
\newcommand{\PGL}{{ PGL}}
\newcommand{\MT}{{ MT}}
\newcommand{\Hdg}{{  Hdg}}
\newcommand{\MTV}{{  MTV}}
\newcommand{\SO}{{  SO}}
\newcommand{\Sp}{{ Sp}}
\newcommand{\Hom}{{ Hom}}
\newcommand{\Ker}{{ Ker}}
\newcommand{\Lie}{{ Lie}}
\newcommand{\Aut}{{ Aut}}
\newcommand{\Image}{{ Image}}
\newcommand{\Gr}{{ Gr}}
\newcommand{\Id}{{Id}}
\newcommand{\rk}{{ rk}}
\newcommand{\pardeg}{{par.deg}}
\newcommand{\SU}{{ SU}}
\newcommand{\Res}{{ Res}}
\newcommand{\Fr}{{ Frob_p}}
\newcommand{\Spec}{{ Spec}}
\newcommand{\Ext}{{ Ext}}
\newcommand{\Sym}{{ Sym}}
\newcommand{\Tor}{{Tor}}
\newcommand{\ch}{{ ch}}
\newcommand{\qtq}[1]{\quad\mbox{#1}\quad}
\newcommand{\bp}{\bar{\partial}}
\newcommand{\Om}{\Omega}
\newcommand{\td}{ Td}
\newcommand{\ind}{Ind}
\newcommand{\ds}{\oplus}
\newcommand{\bds}{\bigoplus}
\newcommand{\ts}{\otimes}
\newcommand{\bts}{\bigotimes}
\newcommand{\diag}{ diag}
\newcommand{\st}{\stackrel}
\newcommand{\btheorem}{\begin{theorem}}
\newcommand{\etheorem}{\end{theorem}}
\newcommand{\bproposition}{\begin{proposition}}
\newcommand{\eproposition}{\end{proposition}}
\newcommand{\bdefinition}{\begin{definition}}
\newcommand{\edefinition}{\end{definition}}
\newcommand{\bcorollary}{\begin{corollary}}
\newcommand{\ecorollary}{\end{corollary}}
\newcommand{\bproof}{\begin{proof}}
\newcommand{\eproof}{\end{proof}}
\newcommand{\bremark}{\begin{remark}}
\newcommand{\eremark}{\end{remark}}
\newcommand{\eexample}{\end{example}}
\newcommand{\bexample}{\begin{example}}
\newcommand{\la}{\langle}
\newcommand{\elemma}{\end{lemma}}
\newcommand{\blemma}{\begin{lemma}}
\newcommand{\ra}{\rangle}
\newcommand{\sq}{\sqrt{-1}}
\newcommand{\suml}{\sum\limits}
\newcommand{\jk}{dz_{j}\wedge d\bar{z}_{k}}
\newcommand{\ckm}{compact K\"{a}hler manifold\ }
\newcommand{\km}{K\"{a}hler manifold\ }
\newcommand{\hvb}{holomorphic vector bundle\ }
\newcommand{\hhvb}{hermitian holomorphic vector bundle\ }
\newcommand{\hlb}{holomorphic line bundle\ }
\newcommand{\hm}{hermitian manifold\ }
\newcommand{\hpd}{hermitian positive definite\ }
\newcommand{\KM}{K\"{a}hler metric\ }
\newcommand{\ob}{orthonormal basis\ }
\newcommand{\ccm}{compact complex manifold\ }
\newcommand{\hhlb}{hermitian holomorphic line bundle\ }
\newcommand{\chm}{compact hermitian manifold\ }
\newcommand{\ep}{\epsilon}
\newcommand{\om}{\omega}
\newcommand{\Llr}{\Longleftrightarrow}
\newcommand{\Lr}{\Longrightarrow}
\newcommand{\ldo}{linear differential operator\ }
\newcommand{\Supp}{ Supp}
\newcommand{\dsl}{\bigoplus\limits}
\newcommand{\wt}{\widetilde}

\newcommand{\Psh}{ Psh}
\newcommand{\sumo}{\widetilde{\suml}}
\newcommand{\grad}{ grad}
\newcommand{\emb}{\hookrightarrow}
\newcommand{\ut}{\la A_{E,\omega}u, u\ra}
\newcommand{\gt}{\la A^{-1}_{E,\omega}g, g\ra}
\newcommand{\pq}{\wedge^{p,q}T^{*}X\ts E}
\newcommand{\nq}{\wedge^{n,q}T^{*}X\ts E}
\newcommand{\pn}{\wedge^{p,n}T^{*}X\ts E}

\newcommand{\p}{\partial}
\newcommand{\pic}{ Pic}
\newcommand{\ord}{ord}
\newcommand{\Div}{ Div}

\renewcommand{\bar}{\overline}
\newcommand{\eps}{\varepsilon}
\newcommand{\pa}{\partial}
\renewcommand{\phi}{\varphi}
\newcommand{\OO}{\mathcal O}

\newcommand{\ka}{K\"ahler }
\newcommand{\kar}{K\"ahler-Ricci}
\newcommand{\ee}{\end{eqnarray*}}
\newcommand{\be}{\begin{eqnarray*}}
\newcommand{\red}{\textcolor[rgb]{1.00,0.00,0.00}}
\newcommand{\yel}{\textcolor[rgb]{0.00,0.00,1.00}}

\newcommand{\rt}{\right}
\newcommand{\lt}{\left}

\newcommand{\aver}[1]{-\hskip-0.35cm\int_{#1}}
\newcommand{\beq}{\begin{equation}}
\newcommand{\eeq}{\end{equation}}

\newcommand{\bd}{\begin{enumerate}}
\newcommand{\ed}{\end{enumerate}}
\newcommand{\ul}{\underline}
\renewcommand{\hat}{\widehat}
\renewcommand{\tilde}{\widetilde}

\newcommand{\rw}{\rightarrow}
\renewcommand{\bf}{\textbf}
\renewcommand{\sc}{\textsc}
\renewcommand{\it}{\textit}
\newcommand{\md}{\textmd}
\renewcommand{\sf}{\textsf}
\renewcommand{\rm}{\textrm}
\newcommand{\lrw}{\Longrightarrow}

\newcommand{\var}{\varnothing}

\newcommand{\dra}{\dashrightarrow}
\renewcommand{\bf}{\textbf}

\newcommand{\TO}{\Longrightarrow}
\newcommand{\OT}{\Longleftarrow}
\newcommand{\col}{\textcolor[rgb]{0.80,0.10,0.40}}
%%%%%%%%%%%%%%%%%%%%%%%%%%%%%%%%%%%%%%%%%%%%%%%%%%%%%%%%%%%
\renewcommand{\Hom}{\text{Hom}}

\renewcommand{\>}{\rightarrow}

\def\baselinestretch{1}

\title{Positivity and vanishing theorems for ample vector bundles }
\author{Kefeng Liu$^{\dagger}$, Xiaofeng Sun$^{*}$ and Xiaokui Yang$^{\dagger}$}

\date{}
\maketitle
%-Title
%\begin{abstract}

%\end{abstract}
%-Abstract

\begin{abstract}

In this paper, we study the Nakano-positivity and
dual-Nakano-positivity of certain adjoint vector bundles associated
to  ample vector bundles. As applications, we get new vanishing
theorems about ample vector bundles. For example, we prove that if
$E$ is an ample vector bundle  over a compact K\"ahler manifold $X$,
$S^kE\ts \det E$ is both Nakano-positive and dual-Nakano-positive
for any $k\geq 0$. Moreover, $H^{n,q}(X,S^kE\ts \det
E)=H^{q,n}(X,S^kE\ts \det E)=0$ for any $q\geq 1$. In particular, if
$(E,h)$ is a Griffiths-positive vector bundle, the naturally induced
Hermitian vector bundle $(S^kE\ts \det E, S^kh\ts \det h)$ is both
Nakano-positive and dual-Nakano-positive for any $k\geq 0$.

%; $S^kE$ is Nakano-positive and dual-Nakano-positive for
%large $k$. In particular, if $(E,h)$ is a Griffiths-positive vector
%bundle, then $(S^kE,S^kh)$ is Nakano-positive and
%dual-Nakano-positive for large $k$. For $\P^n$, we show that
%$(S^kT\P^n, S^kh_{FS})$ is Nakano-positive and dual-Nakano-positive
%for any $k\geq 2$ where $h_{FS}$ is the standard Fubini-Study
%metric.

\end{abstract}

\section{Introduction}
Let $E$ be a holomorphic vector bundle with a Hermitian metric $h$.
Nakano in \cite{Na} introduced
 an analytic notion of positivity by using the curvature of $(E,h)$, and now it is called
 Nakano positivity. Griffiths defined in \cite{G}   Griffiths positivity of $(E,h)$.
 On a Hermitian line bundle, these two concepts are the same. In general,
 Griffiths positivity is weaker than Nakano positivity. On the other hand, Hartshorne  defined in \cite{Har} the ampleness of a
 vector bundle over a projective manifold.  A vector bundle $E$ is said to be ample if the tautological
 line bundle $\sO_{\P(E^*)}(1)$ is ample over $\P(E^*)$.
 For a line bundle, it is well-known that the ampleness of the bundle
is equivalent to its Griffiths positivity.
   In \cite{G}, Griffiths conjectured that this equivalence is also
valid for vector bundles, i.e. $E$ is an ample vector bundle
   if and only if $E$ carries a Griffiths-positive metric.  As is well-known if $E$ admits a Griffiths-positive metric,
   then $\sO_{\P(E^*)}(1)$ has a Griffiths-positive metric(see Proposition \ref{griffithsample}).
   Finding a Griffiths-positive metric on an ample vector bundle seems
    to be very difficult but  is worth being investigated. In \cite{CF}, Campana and Flenner  gave an affirmative answer
to the Griffiths conjecture when the base $S$ is a projective curve,
see also \cite{U}.
   In \cite{SY}, Siu and Yau proved the Frankel conjecture that every compact K\"ahler manifold with positive holomorphic
   bisectional
   curvature is biholomorphic to the projective space. The positivity of holomorphic bisectional curvature
   is the same as Griffiths positivity of the holomorphic
   tangent bundle. On the other hand, S. Mori(\cite{Mo}) proved  the Hartshorne conjecture that any algebraic manifold with ample tangent vector bundle is
biholomorphic to the projective space.

    In this paper, we consider the existence of positive metrics on
    ample vector bundles. It is well-known that metrics with good curvature
properties are  bridges
    between complex algebraic geometry and complex analytic
    geometry. Various vanishing theorems about ample vector bundles can be found in \cite{D1}, \cite{PLS}, \cite{M},
\cite{SS}, \cite{LN2} and   \cite{LN}. In this paper we take a
different approach, we will construct Nakano-positive and
dual-Nakano-positive metrics on various vector bundles associated to
ample vector bundles.

Let $E$ be a holomorphic vector bundle over a compact K\"ahler
manifold $S$ and $F$ a line bundle over $S$. Let $r$ be the rank of
$E$ and $n$ be the complex dimension of $S$. In the following we
briefly describe our main results.

 \btheorem\label{main0}For any integer $k\geq 0$, if $S^{r+k}E\ts
\det E^*\ts F $ is ample over $S$, then $S^kE\ts F$ is both
Nakano-positive and dual-Nakano-positive. \etheorem

\noindent Here we make no assumption on $E$ and we allow $E$ to be
negative. For definitions about Nakano-positivity,
dual-Nakano-positivity and ampleness,  see Section $\bf{2}$. As
pointed out by Berndtsson
 the Nakano positive part of Theorem \ref{main0}  is a special case of \cite{Bo1} where he proves it in the case of a general
holomorphic fibration, but his method can not give the
dual-Nakano-positive part of Theorem \ref{main0}.  Note that
Nakano-positive vector bundles are not necessarily
dual-Nakano-positive and vice versa. For example, for any $n\geq 2$,
the Fubini-Study metric $h_{FS}$ on the holomorphic tangent bundle
$T\P^n$ of $\P^n$ is semi-Nakano-positive and dual-Nakano-positive.
It is well-known that $T\P^n$ does not admit a smooth Hermitian
metric with Nakano-positive curvature for any $n\geq 2$. It is also
easy to see that the holomorphic {cotangent bundle} of a complex
hyperbolic space form is Nakano-positive and is not
dual-Nakano-positive. On the other hand, by the dual
Nakano-positivity, we can get various new vanishing theorems of type
$H^{q,n}$, see Theorem \ref{vanishing} and Proposition
\ref{induction}.

 As applications of Theorem \ref{main0}, we get the following
results:

 \btheorem\label{main2} Let $E$ be an ample vector bundle
over $S$. \bd\item If $F$ is a  {nef line bundle}, then there exists
$k_0=k_0(S,E)$ such that $S^{k}E\ts F$ is Nakano-positive and
dual-Nakano-positive for any $k\geq k_0$. In particular, $S^kE$ is
Nakano-positive and dual-Nakano-positive for any $k\geq k_0$.

\item If $F$ is an arbitrary {vector bundle}, then there
exists $k_0=k_0(S,E,F)$ such that for any $k\geq k_0$, $S^{k}E\ts F$
is Nakano-positive and dual-Nakano-positive. \ed \noindent Moreover,
if the Hermitian vector bundle $(E,h)$ is  Griffiths-positive, then
for large $k$, $(S^kE, S^k h)$ is Nakano-positive and
dual-Nakano-positive. \etheorem

\noindent The following results follow immediately from Theorem
\ref{main0} and Theorem \ref{main2}:

 \bcorollary\label{3} Let $E$ be
a holomorphic vector  bundle over $S$.

\bd\item If $E$ is ample, $S^{k}E\ts \det E$ is Nakano-positive and
dual-Nakano-positive for any $k\geq 0$.

\item If $E$ is ample and its rank  $r$ is greater than $1$, then $S^mE^*\ts (\det
E)^t$ is Nakano-positive and dual-Nakano-positive  for any $t\geq
r+m-1$.

\item If $S^{r+1}E\ts \det E^*$ is  ample, then $E$ is
Nakano-positive and dual-Nakano-positive, so it is
Griffiths-positive. \ed \ecorollary

\noindent If $(E,h)$ is a Griffiths-positive vector bundle,
Demailly-Skoda proved that $E\ts \det E$ and $E^*\ts (\det E)^r$ are
Nakano-positive if $r>1$(\cite{DS}). Recently, Berndtsson proved in
\cite{Bo1} that $S^kE\ts \det E$ is Nakano-positive as soon as $E$
is ample. For more related results, we refer the reader to
\cite{Bo1}, \cite{Bo2} \cite{Bo3}, \cite{MT1}, \cite{MT2} and
\cite{S} and references therein.

 Let $h_{FS}$ be the Fubini-Study metric on $T\P^n$ and
$S^k h_{FS}$  the induced metric on $S^kT\P^n$ by Veronese mapping.
Let $n\geq 2$. It is easy to see that $T\P^n$ does not admit a
Nakano-positive metric. In particular $(T\P^n, h_{FS})$ is not
Nakano-positive. However, $(S^k T\P^n, S^kh_{FS})$ is
Nakano-positive and dual-Nakano-positive for any $k\geq 2$ since
$(S^{k+n}T\P^n\ts K_{\P^n}, S^{k}h_{FS}\ts \det(h_{FS})^{-1})$ is
Griffiths-positive. This  can be viewed as an  evidence of
positivity of some adjoint vector bundles, namely, vector bundles of
type $S^kE\ts (\det E)^\ell\ts K_S$.

\btheorem\label{a} Let $E$ be an ample vector bundle over $S$. Let
$r$ be the rank of $E$ and $n$ the dimension of $S$.  If $r>1$, then
\bd\item $S^{k}E\ts (\det E)^2\ts K_S$ is Nakano-positive and
dual-Nakano-positive for any $k\geq \max\{n-r, 0\}$. Moreover, the
lower bound is sharp.
\item
$E\ts (\det E)^k\ts K_S$ is Nakano-positive and dual-Nakano-positive
for any $k\geq \max\{n+1-r, 2\}$. Moreover, the lower bound is
sharp. \ed \etheorem

\noindent In general, $\det E\ts K_S$ is not an ample line bundle,
for example, $(S,E)=(\P^3, \sO_{\P^3}(1)\ds \sO_{\P^3}(1))$.
Similarly, in the case $n+1-r>2$, i.e. $1<r<n-1$, the vector bundle
$K_S\ts \left(\det E\right)^{n-r}$ can be a negative line bundle,
for example $(S,E)=(\P^4, \sO_{\P^4}(1)\ds\sO_{\P^4}(1))$. So
Theorem \ref{a} is independent of the (dual-)Nakano-positivity of
$S^{k}E\ts \det E$.

Vanishing theorems follow immediately from the Nakano-positive and
dual-Nakano-positive metrics in Theorem \ref{main0}, Theorem
\ref{main2}, Corollary \ref{3} and Theorem \ref{a}. We discuss them
in Theorem \ref{vanishing}, Proposition \ref{induction} and
Corollary \ref{arbitrary}. In the following, we only state one for
example.

 \bproposition  If $E$ is ample over a compact K\"ahler manifold $X$,
$$H^{n,q}(X,S^kE\ts \det E)=H^{q,n}(X,S^kE\ts \det E)=0$$ for any $q\geq
1$ and $k\geq 0$. \eproposition \noindent It is a generalization of
 Griffiths' vanishing theorem (\cite{G}, Theorem G).

\section{Background material}

Let $E$ be a holomorphic vector bundle over a compact K\"ahler
manifold $S$ and $h$ a Hermitian metric on $E$. There exists a
unique connection $\nabla$ which is compatible with the
 metric $h$ and complex structure on $E$. It is called the Chern connection of $(E,h)$. Let $\{z^i\}_{i=1}^n$ be  local holomorphic coordinates
  on $S$ and  $\{e_\alpha\}_{\alpha=1}^r$ be a local frame
 of $E$. The curvature tensor $R^\nabla\in \Gamma(S,\Lambda^2T^*S\ts E^*\ts E)$ has the form
 \beq R^\nabla=\frac{\sq}{2\pi} R_{i\bar j\alpha}^\gamma dz^i\wedge d\bar z^j\ts e^\alpha\ts e_\gamma\eeq
where $R_{i\bar j\alpha}^\gamma=h^{\gamma\bar\beta}R_{i\bar
j\alpha\bar \beta}$ and \beq R_{i\bar j\alpha\bar\beta}= -\frac{\p^2
h_{\alpha\bar \beta}}{\p z^i\p\bar z^j}+h^{\gamma\bar
\delta}\frac{\p h_{\alpha \bar \delta}}{\p z^i}\frac{\p
h_{\gamma\bar\beta}}{\p \bar z^j}\eeq Here and henceforth we
sometimes adopt the Einstein convention for summation.

\bdefinition
 A Hermitian vector bundle
$(E,h)$ is said to be \emph{Griffiths-positive}, if for any nonzero
vectors $u=u^i\frac{\p}{\p z^i}$ and $v=v^\alpha e_\alpha$,  \beq
\sum_{i,j,\alpha,\beta}R_{i\bar j\alpha\bar \beta}u^i\bar u^j
v^\alpha\bar v^\beta>0\eeq $(E,h)$ is said to be
\emph{Nakano-positive}, if for any nonzero vector
$u=u^{i\alpha}\frac{\p}{\p z^i}\ts e_\alpha$, \beq
\sum_{i,j,\alpha,\beta}R_{i\bar j\alpha\bar \beta} u^{i\alpha}\bar
u^{j\beta}>0 \eeq $(E,h)$ is said to be \emph{dual-Nakano-positive},
if for any nonzero vector $u=u^{i\alpha}\frac{\p}{\p z^i}\ts
e_\alpha$, \beq \sum_{i,j,\alpha,\beta}R_{i\bar j\alpha\bar \beta}
u^{i\beta}\bar u^{j\alpha}>0 \eeq It is easy to see that $(E,h)$ is
dual-Nakano-positive if and only if $(E^*,h^*)$ is Nakano-negative.
\noindent The notions of semi-positivity, negativity and
semi-negativity can be defined similarly. We say $E$ is
Nakano-positive (resp. Griffiths-positive, dual-Nakano-positive,
$\cdots$), if it admits a Nakano-positive(resp. Griffiths-positive,
dual-Nakano-positive, $\cdots$) metric. \edefinition

The following geometric definition of nefness is due to \cite{DPS}.

\bdefinition
 Let $(S,\omega_0)$ be a compact K\"ahler manifold. A line bundle $L$ over $S$ is said to be nef,
if  for any $\eps>0$, there exists a smooth Hermitian metric
$h_\eps$ on $L$ such that the curvature of $(L,h_\eps)$ satisfies
\beq R=-\frac{\sq}{2\pi}\p\bp\log h_\eps\geq -\eps\omega_0\eeq
\edefinition

 \noindent This means that the curvature of $L$ can have an arbitrarily small negative
part. Clearly a nef line bundle $L$ satisfies $$\int_C c_1(L)\geq
0$$ for any irreducible curve $C\subset S$.
 For projective
algebraic $S$, both notions coincide.

By the Kodaira embedding theorem, we have the following geometric
definition of ampleness.

\bdefinition Let $(S,\omega_0)$ be a compact K\"ahler manifold. A
line bundle $L$ over $S$ is said to be ample, if there exists a
smooth Hermitian metric $h$ on $L$ such that the curvature $R$ of
$(L,h)$ satisfies \beq R=-\frac{\sq}{2\pi}\p\bp\log h>0\
 \eeq

\edefinition

 \noindent For  comprehensive
descriptions of positivity, nefness, ampleness and related topics,
see \cite{D}, \cite{DPS}, \cite{L}, \cite{G}, \cite{SS} and
\cite{U}.

 Let $E$ be a Hermitian vector bundle of rank $r$ over a compact K\"ahler manifold $S$, $L=\sO_{\P(E^*)}(1)$ be the
  tautological line bundle of the projective bundle $\P(E^*)$ and $\pi$  the canonical projection $\P(E^*)\> S$.
  By definition(\cite{Har}), $E$ is an ample vector bundle over $S$ if
   $\sO_{\P(E^*)}(1)$ is an ample line bundle over $\P(E^*)$. $E$ is said to be
   nef, if $\sO_{\P(E^*)}(1)$ is nef.
   To simplify the
   notations we will denote $\P(E^*)$ by $X$ and the fiber
   $\pi^{-1}(\{s\})$ by $X_s$.

Let $(e_1,\cdots, e_r)$ be the local holomorphic frame with respect
to a given trivialization on $E$ and the dual frame on $E^*$ is
denoted by $(e^1,\cdots, e^r)$. The corresponding holomorphic
coordinates on $E^*$ are denoted by $(W_1,\cdots, W_r)$.  There is a
local section $e_{L^*}$ of $L^*$ defined by \beq
e_{L^*}=\sum_{\alpha=1}^r W_\alpha e^\alpha\eeq Its dual section is
denoted by $e_L$. Let $h^E$ be a fixed Hermitian metric on $E$ and
$h^L$ the induced quotient metric by the morphism
$(\pi^*E,\pi^*h^E)\>L$. %

If $\left(h_{\alpha\bar\beta}\right)$ is the matrix representation
of $h^E$ with respect to the basis $\{e_\alpha\}_{\alpha=1}^r$, then
$h^L$ can be written as \beq
h^L=\frac{1}{h^{L^*}(e_{L^*},e_{L^*})}=\frac{1}{\sum
h^{\alpha\bar\beta}W_\alpha\bar W_\beta}\label{inducedmetric} \eeq

\bproposition\label{quotient} The curvature of $(L,h^L)$ is \beq
R^{h^L}=-\frac{\sq}{2\pi}\p\bp\log
h^L=\frac{\sq}{2\pi}\p\bp\log\left(\sum
h^{\alpha\bar\beta}W_\alpha\bar W_\beta\right)
\label{inducedcurvature}\eeq where $\p$ and $\bp$ are operators on
the total space $\P(E^*)$. \eproposition

Although the following result is well-known(\cite{D}, \cite{G}), we
include a proof here for the sake of completeness.
\bproposition\label{griffithsample} If $(E,h^E)$ is a
Griffiths-positive vector bundle, then $E$ is ample. \bproof  We
will show that
 the induced metric  $h^L$  in (\ref{inducedmetric})
is positive. We fix a point $p\in \P(E^*)$, then there exist local
holomorphic coordinates
 $(z^1,\cdots, z^n)$ centered at point $s=\pi(p)$ and local holomorphic basis $\{e_1,\cdots, e_r\}$ of $E$ around $s$ such that
 \beq h_{\alpha\bar\beta}=\delta_{\alpha\bar\beta}-R_{i\bar j \alpha\bar\beta}z^i\bar z^j+O(|z|^3) \label{normal}\eeq
Without loss of generality, we assume $p$  is the point $(0,\cdots,
0,[a_1,\cdots, a_r])$ with $a_r=1$. On the chart $U=\{W_r=1\}$ of
the fiber $\P^{r-1}$, we set $w^A=W_A$ for $A=1,\cdots, r-1$. By
formula (\ref{inducedcurvature}) and (\ref{normal}) \beq
R^{h^L}(p)=\frac{\sq}{2\pi}\left(\sum R_{i\bar
j\alpha\bar\beta}\frac{a_\beta \bar a_\alpha}{|a|^2}dz^i\wedge d\bar
z^j+\sum_{A,B=1}^{r-1}\left(1-\frac{a_B\bar
a_A}{|a|^2}\right)dw^A\wedge d\bar w^B\right) \eeq where
$|a|^2=\sum\limits_{\alpha=1}^r|a_\alpha|^2$. If $R^{E}$ is Griffith
positive,$$\left(\sum_{\alpha,\beta=1}^rR_{i\bar
j\alpha\bar\beta}\frac{a_\beta \bar a_\alpha}{|a|^2}\right)$$ is a
Hermitian positive $n\times n$ matrix. Consequently, $R^{h^L}(p)$ is
a Hermitian positive $(1,1)$ form on $\P(E^*)$, i.e. $h^L$ is a
positive Hermitian metric. \eproof \eproposition

 The following linear algebraic lemma will be used in Theorem \ref{Nakano}.
\blemma\label{linear} If the matrix

$$
T=\left(%
\begin{array}{cc}
 A &B \\
  C & D \\
\end{array}%
\right) $$ is invertible and $D$ is invertible, then
$(A-BD^{-1}C)^{-1}$ exists and
$$
T^{-1}=\left(%
\begin{array}{cc}
  (A-BD^{-1}C)^{-1}&-(A-BD^{-1}C)^{-1}BD^{-1} \\
  -D^{-1}C(A-BD^{-1}C)^{-1} & D^{-1}C(A-BD^{-1}C)^{-1}BD^{-1}+D^{-1}\\
\end{array}%
\right) $$ Moreover, if $T$ is positive definite, then $A-BD^{-1}C$
is positive definite.\elemma

\section{Curvature formulas}

Let $F$ be a holomorphic line bundle over $S$, $L=\sO_{\P(E^*)}(1)$
and $\pi:\P(E^*)\>S$. For simplicity of notations,  we set $\tilde
L=L^k\ts \pi^*( F)$ for $k\geq 0$ and $X=\P(E^*)$. Let $h_0$ be a
Hermitian metric
    on $\tilde L$ and $\{\omega_s\}_{s\in S}$ a smooth family of K\"ahler metrics
    on the fibers $X_s=\P(E_s^*)$ of $X$ which are induced by the curvature form of some metric on $\sO_{\P(E^*)}(1)$.  Let $\{w^A\}_{A=1}^{r-1}$
    be the local holomorphic coordinates on the fiber $X_s$ which are induced by the homogeneous coordinates $[W_1,\cdots,
    W_r]$ on a trivialization chart. Using these notations, we can write  $\omega_s$
    as
\beq \omega_s=\frac{\sq }{2\pi}\sum_{A,B=1}^{r-1} g_{A\bar B}(s,w)
dw^A\wedge d\bar w^B \label{ws}\eeq
  It is well-known $H^0(\P^{r-1}, \sO_{\P^{r-1}}(k))$ can be identified as the space of
 homogeneous polynomials of degree $k$ in $r$ variables. Therefore, the sections of $H^0(X_s, \tilde L|_{X_s})$ are of the form $V_\alpha e_{L}^{\ts k}\ts \underline e$
 where $V_\alpha$ are homogenous polynomials in $\{W_1,\cdots ,W_r\}$ of degree $k$ and $\underline e$ the base of $\pi^*(F)$ induced by a base $e$ of $F$.
 For example, if $\alpha=(\alpha_1,\cdots,\alpha_r)$ with $\alpha_1+\cdots+\alpha_r=k$ and $\alpha_j$  are nonnegative integers,
 \beq V_\alpha=W_1^{\alpha_1}\cdots W_r^{\alpha_r}.\label{valpha}\eeq

\noindent  Now we set
  $$E_\alpha=e_1^{\ts\alpha_1}\ts\cdots\ts e_r^{\ts\alpha_r}\ts e\qtq{and} e_{\tilde L}=e_L^{\ts k}\ts\underline e$$
  which are bases of $S^kE\ts F$ and $\tilde L$ respectively.
  We obtain a vector bundle whose fibers are $H^0(X_s, \tilde L|_{X_s})$. In fact, this vector bundle is $\tilde E=S^kE\ts F$.
  Now we can define a smooth Hermitian metric $f$ on $S^kE\ts F$  by $(\tilde L,h_0)$ and $(
 X_s,\omega_s)$, locally it is
    \begin{eqnarray} f_{\alpha\bar\beta}:=f(E_\alpha,E_\beta)\nonumber&=&\int_{X_s} \langle V_\alpha e_{\tilde L}, V_\beta e_{\tilde L}\rangle_{h_0}\frac{\omega_s^{r-1}}{(r-1)!}
    \\&=&\int_{X_s} h_0 V_\alpha\bar V_\beta \frac{\omega_s^{r-1}}{(r-1)!}\label{l2metric}
    \end{eqnarray}
Here we regard $h_0$ locally as a positive function.  In this
general setting, the Hermitian metric $h_0$ on $\tilde L$ and
K\"ahler metrics $\omega_s$ on the fibers are independent.

Let $(z^1,\cdots, z^n)$ be  local holomorphic coordinates on $S$. By
definition, the curvature tensor of $f$  is \beq R_{i\bar
j\alpha\bar\beta}=-\frac{\p^2 f_{\alpha\bar \beta}}{\p z^i\p\bar
z^j}+\sum_{\gamma,\delta}f^{\gamma\bar \delta}\frac{\p f_{\alpha
\bar \delta}}{\p z^i}\frac{\p f_{\gamma\bar\beta}}{\p \bar
z^j}\label{curvaturedefinition}\eeq  In the following, we will
compute the curvature of $f$. Let $T_{X/S}$ be the relative tangent
bundle of the fibration $\P(E^*)\rightarrow S$, then $g_{A\bar B}$
is a metric on $T_{X/S}$ and $\det (g_{A\bar B})$ is a metric on
$\det(T_{X/S})$.
 Let
$\phi=-\log(h_0\det(g_{A\bar B}))$ be the local weight of induced
Hermitian metric $h_0\det(g_{A\bar B})$ on $\tilde L\ts
\det(T_{X/S})$. In the sequel, we will use the following notations
$$\phi_i=\frac{\p\phi}{\p z^i},  \phi_{i\bar j}=\frac{\p^2\phi}{\p z^i\p\bar z^j}, \phi_{A\bar B}=\frac{\p^2 \phi}{\p w^A \p \bar w^B},
\phi_{i\bar B}=\frac{\p^2 \phi}{\p z^i\p \bar w^B}, \phi_{A\bar
j}=\frac{\p^2\phi}{\p \bar z^j \p w^A}$$ and $(\phi^{A\bar B})$ is
the transpose inverse of the $(r-1)\times (r-1)$ matrix
$(\phi_{A\bar B})$,
$$\suml_{B=1}^{r-1}\phi^{A\bar B}\phi_{C\bar B}=\delta_C^A$$

The following lemma can be deduced from the formulas in \cite{S1},
\cite{Wo} and \cite{Siu}. In the case of holomorphic fibration $\P(
E^*)\rightarrow S$, we can compute it directly.
 \blemma\label{firstderivative} The first order derivative of $f_{\alpha\bar\beta}$ is
\beq \frac{\p f_{\alpha\bar \beta}}{\p z^i}=-\int_{X_s} h_0 V_\alpha
\bar V_\beta \phi_i \frac{\omega_s^{r-1}}{(r-1)!}=\int_{X_s}\langle
-V_\alpha\phi_ie_{\tilde L}, V_\beta e_{\tilde
L}\rangle_{h_0}\frac{\omega_s^{r-1}}{(r-1)!} \eeq \bproof By the
local expression (\ref{ws}) of $\omega_s$,
$$\frac{\omega_s^{r-1}}{(r-1)!}=\det(g_{A\bar B}) dV_{\C^{r-1}}$$
where $dV_{\C^{r-1}}$ is standard volume on $\C^{r-1}$. Therefore
$$f_{\alpha\bar\beta}=\int_{X_s}e^{-\phi} V_\alpha\bar V_\beta dV_{\C^{r-1}}$$
and the first order derivative is \be \frac{\p f_{\alpha\bar
\beta}}{\p z^i}&=&\int_{X_s}\frac{\p
e^{-\phi}}{\p z^i} V_\alpha\bar V_\beta dV_{\C^{r-1}}\\
&=& -\int_{X_s} \phi_i e^{-\phi} V_\alpha\bar V_\beta
dV_{\C^{r-1}}\\&=&-\int_{X_s} h_0 V_\alpha \bar V_\beta \phi_i
\frac{\omega_s^{r-1}}{(r-1)!}\ee
 \eproof \elemma

\btheorem\label{curvature} The curvature tensor of the Hermitian
metric $f$ on $S^kE\ts F$ is  \beq R_{i\bar j \alpha\bar
\beta}=\int_{ X_s} h_0 V_\alpha\bar V_\beta \phi_{i\bar
j}\frac{\omega_s^{r-1}}{(r-1)!} -\int_{X_s} h_0P_{i\alpha}\bar
P_{j\beta}\frac{\omega_s^{r-1}}{(r-1)!} \label{curvatureformula}\eeq
where \beq P_{i\alpha}=-V_\alpha
\phi_i-\sum_{\gamma}^{}V_\gamma\left(\sum_{\delta}f^{\gamma \bar
\delta}\frac{\p f_{\alpha\bar \delta}}{\p z^i}\right) \eeq

\bproof The idea we use is due to Berndtsson(\cite{Bo1}, Section
$2$). For simplicity of notations, we set $A_{i\alpha}=-V_\alpha
\phi_i$. The Hermitian metric (\ref{l2metric}) is also a norm on the
smooth section space $\Gamma(X_s, \tilde L|_{X_s})$, and  it induces
an orthogonal projection
$$\tilde \pi_s: \Gamma(X_s, \tilde L|_{X_s})\> H^0(X_s,\tilde L|_{X_s})$$
\noindent Using this projection, we can rewrite the first order
derivative as  \be \frac{\p f_{\alpha\bar \beta}}{\p z^i}
&=&\int_{X_s}\langle A_{i\alpha}e_{\tilde L}, V_\beta e_{\tilde L}\rangle_{h_0}\frac{\omega_s^{r-1}}{(r-1)!}\\
&=&\int_{X_s}\langle \tilde \pi_s(A_{i\alpha}e_{\tilde L})+(A_{i\alpha}e_{\tilde L}-\tilde\pi_s(A_{i\alpha}e_{\tilde L})), V_\beta e_{\tilde L}\rangle_{h_0}\frac{\omega_s^{r-1}}{(r-1)!}\\
&=&\int_{X_s}\langle \tilde\pi_s(A_{i\alpha}e_{\tilde L}), V_\beta
e_{\tilde L}\rangle_{h_0}\frac{\omega_s^{r-1}}{(r-1)!}\ee since
$(A_{i\alpha}e_{\tilde L}-\tilde\pi_s(A_{i\alpha}e_{\tilde L}))$ is
in the orthogonal complement of $H^0(X_s,{\tilde L}|_{X_s})$. By
this relation, we can  write $\tilde\pi_s(A_{i\alpha}e_{\tilde L})$
in the basis  $\{V_\alpha e_{\tilde L}\}$ of $H^0(X_s, {\tilde L}|_{
X_s})$,  \beq \tilde \pi_s(A_{i\alpha}e_{\tilde
L})=\sum_{\gamma}\left(\sum_{\delta}f^{\gamma \bar \delta}\frac{\p
f_{\alpha\bar \delta}}{\p z^i}\right)\left(V_\gamma e_{\tilde
L}\right) \eeq

\noindent From this identity, we obtain  \beq \int_{ X_s}
\big\langle \tilde \pi_s(A_{i\alpha} e_{\tilde L}),
\tilde\pi_s(A_{j\beta}e_{\tilde
L})\big\rangle_{h_0}\frac{\omega_s^{r-1}}{(r-1)!}=\sum_{\gamma,\delta}f^{\gamma\bar
\delta}\frac{\p f_{\alpha \bar \delta}}{\p z^i}\frac{\p
f_{\gamma\bar\beta}}{\p \bar z^j}\eeq Suppose \beq
P_{i\alpha}=A_{i\alpha}-\sum_{\gamma}^{}V_\gamma\left(\sum_{\delta}f^{\gamma
\bar \delta}\frac{\p f_{\alpha\bar \delta}}{\p z^i}\right) \eeq then
$ A_{i\alpha} e_{\tilde L}=\tilde\pi_s(A_{i\alpha} e_{\tilde
L})+P_{i\alpha} e_{\tilde L} $, that is, \beq
\tilde\pi_s(P_{i\alpha} e_{\tilde L} )=0\label{0projection} \eeq

\noindent Similar to Lemma \ref{firstderivative}, we obtain the
second order derivative \be \frac{\p^2 f_{\alpha\bar \beta}}{\p
z^i\p \bar z^j}&=&-\int_{X_s} h_0 V_\alpha\bar V_\beta \phi_{i\bar
j}\frac{\omega_s^{r-1}}{(r-1)!} +\int_{X_s} \langle
V_\alpha \phi_i e_{\tilde L}, V_\beta \phi_j e_{\tilde L}\rangle_{h_0}\frac{\omega_s^{r-1}}{(r-1)!}\\
&=&-\int_{X_s} h_0 V_\alpha\bar V_\beta \phi_{i\bar
j}\frac{\omega_s^{r-1}}{(r-1)!}+\int_{X_s} \langle
A_{i\alpha}e_{\tilde L}, A_{j\beta} e_{\tilde L}\rangle_{h_0}\frac{\omega_s^{r-1}}{(r-1)!}\\
&=&-\int_{X_s} h_0 V_\alpha\bar V_\beta \phi_{i\bar
j}\frac{\omega_s^{r-1}}{(r-1)!}\\&&+\int_{X_s} \big\langle
P_{i\alpha}e_{\tilde L} +\tilde \pi_s(A_{i\alpha} e_{\tilde L}), P_{j\beta} e_{\tilde L}+\tilde\pi_s(A_{j\beta}e_{\tilde L})\big\rangle_{h_0}\frac{\omega_s^{r-1}}{(r-1)!}\\
&=&-\int_{X_s} h_0 V_\alpha\bar V_\beta \phi_{i\bar
j}\frac{\omega_s^{r-1}}{(r-1)!} +\int_{X_s} h_0P_{i\alpha}\bar
P_{j\beta}\frac{\omega_s^{r-1}}{(r-1)!}+ \\&&\int_{X_s} \big\langle
\tilde \pi_s(A_{i\alpha}
e_{\tilde L}), \tilde\pi_s(A_{j\beta}e_{\tilde L})\big\rangle_{h_0}\frac{\omega_s^{r-1}}{(r-1)!}\\
&=&-\int_{X_s} h_0 V_\alpha\bar V_\beta \phi_{i\bar
j}\frac{\omega_s^{r-1}}{(r-1)!} +\int_{X_s} h_0P_{i\alpha}\bar
P_{j\beta}\frac{\omega_s^{r-1}}{(r-1)!}\\&&+ f^{\gamma\bar
\delta}\frac{\p f_{\alpha \bar \delta}}{\p z^i}\frac{\p f_{\gamma\bar\beta}}{\p \bar z^j}\\
 \ee
\noindent By formula (\ref{curvaturedefinition}), we get the
curvature formula (\ref{curvatureformula}).
 \eproof \etheorem

\section{Positivity of Hermitian metrics}

 If $(E,h)$ is a Griffiths-positive, Demailly-Skoda(\cite{DS}) showed
 that $(E\ts\det E,h\ts \det h)$ is Nakano-positive. They proved it
 by using a discrete Fourier transformation method. Here, we use a linear algebraic
 argument to show  $(E\ts \det E, h\ts \det h)$ is both Nakano-positive and
 dual-Nakano-positive.

Let $\omega_{FS}$ be the standard Fubini-Study metric on $\P^{r-1}$
and
 $[W_1,\cdots W_r]$ the homogeneous coordinates on $\P^{r-1}$. If $A=(\alpha_1,\cdots, \alpha_k)$
and $B=(\beta_1,\beta_2,\cdots, \beta_k)$, we  define the
generalized Kronecker-$\delta$ for multi-index by the following
formula \beq \delta_{AB}=\sum_{\sigma\in S_k}
\prod_{j=1}^k\delta_{\alpha_{\sigma(j)}\beta_{\sigma(j)}}
\label{generalizeddelta}\eeq where $S_k$ is the permutation group in
$k$ symbols. \blemma\label{calculation} If $V_A=W_{\alpha_1}\cdots
W_{\alpha_k}$and $V_B=W_{\beta_1}\cdots W_{ \beta_k}$, then
\beq\int_{\P^{r-1}} \frac{V_{A}\bar
V_B}{|W|^{2k}}\frac{\omega_{FS}^{r-1}}{(r-1)!}=\frac{\delta_{AB}}{(r+k-1)!}
\eeq For simple-index notations, \beq\int_{\P^{r-1}} \frac{W_\alpha
\bar
W_\beta}{|W|^2}\frac{\omega_{FS}^{r-1}}{(r-1)!}=\frac{\delta_{\alpha\beta}}{r!}
,\  \int_{\P^{r-1}}\frac{W_\alpha \bar{W_\beta} W_\gamma
\bar{W_\delta}}{|W|^4}\frac{\omega_{FS}^{r-1}}{(r-1)!}=\frac{\delta_{\alpha\beta}\delta_{\gamma\delta}+\delta_{\alpha\delta}\delta_{\beta\gamma}}{(r+1)!}
\eeq\elemma

Without loss of generality we can assume, at a fixed $s\in
 S$, $h_{\alpha\bar\beta}(s)=\delta_{\alpha\beta}$. The curvature of
 $(E\ts \det E, h\ts \det h)$ is
 \beq R_{i\bar j\alpha\bar\beta}^{E\ts \det E}(s)=R_{i\bar j \alpha\bar\beta}(s)+\delta_{\alpha\beta}\cdot
\sum_{\gamma}R_{i\bar j \gamma\bar\gamma}(s) \eeq By Lemma
\ref{calculation}, we obtain \beq R_{i\bar j
\alpha\bar\beta}(s)+\delta_{\alpha\beta}\cdot \sum_{\gamma}R_{i\bar
j \gamma\bar\gamma}(s)=r! \cdot \int_{ \P^{r-1}} \frac{ W_\alpha\bar
W_\beta}{|W|^{2}} \phi_{i\bar j}\frac{\omega^{r-1}_{FS}}{(r-1)!}
\label{curv}\eeq where \beq \phi_{i\bar
j}=(r+1)\sum_{\gamma,\delta}R_{i\bar j
\gamma\bar\delta}(s)\frac{W_\delta\bar W_\gamma}{|W|^2} \eeq If
$(E,h)$ is Griffiths-positive, then $(\phi_{i\bar j})$ is Hermitian
positive. For any nonzero $u=(u^{i\alpha})$ \beq R^{E\ts \det
E}_{i\bar j \alpha\bar\beta}u^{i\beta}\bar
u^{j\alpha}=(r+1)\int_{\P^{r-1}}\phi_{i\bar j}\frac{(u^{i\beta}\bar
W_\beta)\cdot \bar{\left(u^{j\alpha}\bar
W_\alpha\right)}}{|W|^2}\frac{\omega^{r-1}_{FS}}{(r-1)!}>0 \eeq
Therefore, $(E\ts\det E, h\ts \det h)$ is dual-Nakano-positive. By a
similar formulation, we know $(E\ts \det E, h\ts \det h)$ is
Nakano-positive. For more related results, see Section $7$.

 In the following, we will prove similar results for ample vector
 bundles.

\subsection{Nakano-positivity}
In this subsection, we will use  $\bp$-estimate on a compact
K\"ahler manifold to analyze the curvature formula in Theorem
\ref{curvature}, $$ R_{i\bar j \alpha\bar \beta}=\int_{\ X_s} h_0
V_\alpha\bar V_\beta \phi_{i\bar j}\frac{\omega_s^{r-1}}{(r-1)!}
-\int_{X_s} h_0P_{i\alpha}\bar
P_{j\beta}\frac{\omega_s^{r-1}}{(r-1)!} $$ The first term on the
right hand side involves the horizontal direction curvature
$\phi_{i\bar j}$ of the line bundle ${\tilde L}\ts \det(T_{X/S})$.
If the line bundle ${\tilde L}\ts \det(T_{X/S})$ is positive in the
horizontal direction, we can choose $(h_0, \omega_s)$ such that
$\phi$ is positive  in the  horizontal direction, i.e. $(\phi_{i\bar
j})$ is Hermitian positive. We will get a  lower bound of the second
term by using H\"ormander's $L^2$-estimate, following an idea of
Berndtsson(\cite{Bo1}).

\blemma\label{L} Let $(M^n,\omega_g)$ be a compact K\"ahler manifold
and $(L,h)$  a Hermitian line bundle over $M$. If there exists a
positive constant $c$ such that \beq Ric(\omega_g)+R^h\geq
c\omega_g\label{condition}\eeq then for any $w\in
\Gamma(M,T^{*0,1}M\ts L)$ such that $\bp w=0$, there exists a unique
$u\in \Gamma(M,L)$ such that $\bp u=w$ and $\tilde\pi(u)=0$ where
$\tilde\pi: \Gamma(M,L)\> H^0(M,L)$ is the orthogonal projection.
Moreover,  \beq \int_M |u|_{h}^2 \frac{\omega_g^n}{n!}\leq
\frac{1}{c}\int_{M}|w|^2_{g^*\ts h}
\frac{\omega_g^n}{n!}\label{estimate}\eeq  \elemma

\noindent We refer the reader to \cite{D} and \cite{H} for the proof
of Lemma \ref{L}.

 Now we apply Lemma \ref{L} to each fiber $( X_s, \omega_s)$
and $({\tilde L}|_{X_s}, h_0|_{X_s})$. At a fixed point $s\in S$,
the fiber direction curvature of the induced metric on $\tilde L\ts
\det(T_{X/S})$ is \beq
-\frac{\sq}{2\pi}\p_s\bp_s\log(h_0\det(g_{A\bar B}))= R^{\tilde
L^{h_0}_s}+Ric_F(\omega_s) \eeq On the other hand
$$-\frac{\sq}{2\pi}\p_s\bp_s\log(h_0\det(g_{A\bar B}))=\frac{\sq}{2\pi} \p_s\bp_s\phi$$
where $\phi=-\log(h_0\det(g_{A\bar B}))$. So condition
(\ref{condition}) turns out to be \beq (\phi_{A\bar B})\geq c_s
(g_{A\bar B}) \label{main1}\eeq for some positive constant
$c_s=c(s)$.

\btheorem\label{b2} If $(\phi_{A\bar B})\geq c_s (g_{A\bar B})$ at
point $s\in S$, then for any
$$u=\sum_{i,\alpha} u^{i\alpha}\frac{\p}{\p z^i}\ts E_\alpha\in \Gamma(S, T^{1,0}S \ts \tilde E)$$
with $\tilde E=S^kE\ts F$, we have the following estimate at point
$s$, \beq R_{i\bar j\alpha\bar \beta} u_{}^{i\alpha}\bar
{u_{}^{j\beta}}\geq \int_{ X_s}h_0( V_\alpha u_{}^{i\alpha})
\bar{(V_\beta u_{}^{j\beta})}\left(\phi_{i\bar j}-\frac{g^{A\bar
B}\phi_{i\bar B}\phi_{A\bar
j}}{c_s}\right)\frac{\omega_s^{r-1}}{(r-1)!}\label{estimate1}\eeq

\bproof At point $s\in S$, we  set
$$P=\sum_{i,\alpha}P_{i\alpha}u^{i\alpha} e_{\tilde L} \in  \Gamma( X_s, {\tilde L}_s),\ \  K= -\sum_{i,\alpha} V_\alpha \phi_i u^{i\alpha} e_{\tilde L} \in \Gamma( X_s, {\tilde L}_s)$$
It is obvious that $\bp_s P=\bp_s K$ where $\bp_s$ is $\bp$ on the
fiber direction. On the other hand,  by (\ref{0projection}),
$\tilde\pi_s(P)=0$. So we can apply Lemma \ref{L} and get \beq
\int_{X_s}|P|^2_{h_0}\frac{\omega_s^{r-1}}{(r-1)!}\leq
\frac{1}{c_s}\int_{X_s}|\bp_s K|_{g^*_s\ts h_0}^2
\frac{\omega_s^{r-1}}{(r-1)!}\label{l2}\eeq Since $\bp_sK=
-\suml_{i,\alpha,B} V_\alpha \phi_{i\bar B} u^{i\alpha} d\bar z^B
\ts e_{\tilde L}$,
$$|\bp_s K|^2_{g_s^*\ts h_0}=\sum_{i,j}\sum_{\alpha,\beta } h_0( V_\alpha u_{}^{i\alpha})
\bar{(V_\beta u_{}^{j\beta})}g^{A\bar B}\phi_{i\bar B}\phi_{A\bar
j}$$  By inequality (\ref{l2}) and Theorem \ref{curvature},  we get
the estimate (\ref{estimate1}).
 \eproof
\etheorem

 Before proving the main theorems, we need the following
lemma:

\blemma\label{positive} If $E$ is a holomorphic vector bundle with
rank $r$ over a compact K\"ahler manifold $S$  and $F$ is a line
bundle over $S$ such that $S^{k+r}E\ts \det E^*\ts F$ is  ample over
$S$, then there exists a positive Hermitian metric $\lambda_0$ on
$\sO_{\P(E^*)}(k)\ts \pi^*(F)\ts \det(T_{X/S})$. \bproof Let $\hat
E$ be $S^{k+r}E\ts \det(E^*)\ts F$.
 It is obvious that
$\P(S^{k+r}E^*)=\P(\hat E^*)$. The tautological line bundles of them
are related by the following formula
 \beq \sO_{\P(\hat E^*)}(1)=\sO_{\P(S^{k+r}E^*)}(1)\ts
\pi^*_{k+r}(\det E^*)\ts\pi^*_{k+r}(F) \eeq where
$\pi_{k+r}:\P(S^{k+r}E^*)\rightarrow S$ is the canonical projection.
Let $v_{k+r}:\P(E^*)\rightarrow\P(S^{k+r}E^*)$ be the standard
Veronese embedding, then  \beq
\sO_{\P(E^*)}(k+r)=v_{k+r}^*\left(\sO_{\P(S^{k+r}E^*)}(1)\right)
\eeq

\noindent Similarly, let $\mu_{k+r}$ be the induced mapping
$\mu_{k+r}:\P(E^*)\> \P(\hat E^*)$, then \beq
\mu_{k+r}^*\left(\sO_{\P(\hat
E^*)}(1)\right)=\sO_{\P(E^*)}(k+r)\ts\pi^*(F\ts\det E^*) \eeq By the
identity \beq K_{ X}=\pi^*(K_S)\ts \sO_{\P(E^*)}(-r)\ts \pi^*(\det
E), \eeq we obtain \beq \mu_{k+r}^*\left(\sO_{\P(\hat
E^*)}(1)\right)=\sO_{\P(E^*)}(k)\ts \pi^*(F)\ts \det(T_{X/S})=\tilde
L\ts \det(T_{X/S})\eeq If $\hat E$ is ample, then $\sO_{\P(\hat
E^*)}(1)$ is ample and so is $\tilde L\ts \det(T_{X/S})$. So there
exists a positive Hermitian metric $\lambda_0$ on $\tilde L\ts
\det(T_{X/S})$. \eproof \elemma

\btheorem\label{Nakano} Let $E$ be a holomorphic vector bundle over
a compact K\"ahler manifold $S$ and $F$  a line bundle over $S$. Let
$r$ be the rank of $E$ and $k\geq 0$ an arbitrary integer. If
$S^{k+r}E\ts \det E^*\ts F$ is ample over $S$, then there exists a
smooth Hermitian metric $f$ on $S^kE\ts F$ such that $(S^kE\ts F,
f)$ is
 Nakano-positive. \bproof By Lemma \ref{positive},
 there exists a positive Hermitian metric $ \lambda_0$   on the ample line bundle $\tilde L\ts \det(T_{X/S})$. We set
$$\omega_s=-\frac{\sq}{2\pi} \p_s\bp_s\log \lambda_0=\frac{\sq}{2\pi}\sum_{A,B=1}^{r-1} g_{A\bar B}(s,w)dw^A\wedge d\bar w^B$$
which is a smooth family of K\"ahler metrics on the fibers $ X_s$.
We get  an induced Hermitian metric on $\tilde L$, namely, \beq
h_0=\frac{\lambda_0}{\det(g_{A\bar B})}\label{45} \eeq

Let $f$ be the Hermitian metric on the vector bundle $S^kE\ts \det
F$ induced by $(\tilde L, h_0)$ and $( X_s, \omega_s)$(see
(\ref{l2metric})). In this setting, the weight $\phi$ of induced
metric on $\tilde L\ts \det(T_{X/S})$ is
$$\phi =-\log\left( h_0 \det( g_{A\bar B})\right)=-\log \lambda_0$$
Hence \beq \left( \phi_{A\bar B}\right)=\left(g_{A\bar B}\right)\eeq
and in Theorem \ref{b2}, $c_s=1$ for any $s\in S$. Therefore \be
R^{\tilde E}(u,u)&=& R_{i\bar j\alpha\bar \beta} u_{}^{i\alpha}\bar
{u_{}^{j\beta}}\\&\geq &\int_{ X_s} h_0( V_\alpha u_{}^{i\alpha})
\bar{(V_\beta u_{}^{j\beta})}\left( \phi_{i\bar
j}-{\suml_{A,B=1}^{r-1}g^{A\bar B} \phi_{i\bar B} \phi_{A\bar
j}}\right)\frac{\omega_s^{r-1}}{(r-1)!}\\&=& \int_{X_s} h_0(
V_\alpha u_{}^{i\alpha}) \bar{(V_\beta u_{}^{j\beta})}\left(
\phi_{i\bar j}-\suml_{A,B=1}^{r-1}\phi^{A\bar B} \phi_{i\bar B}
\phi_{A\bar j}\right)\frac{\omega_s^{r-1}}{(r-1)!}\ee for any
$u=\suml_{i,\alpha} u^{i\alpha}\frac{\p}{\p z^i}\ts E_\alpha\in
\Gamma(S, T^{1,0}S \ts \tilde E)$.

 On the other hand  $\lambda_0$ is a positive Hermitian metric on the line bundle $\tilde L\ts \det(T_{X/S})$.  The curvature form of $\lambda_0$ can be represented by
  a Hermitian positive  matrix, namely, the coefficients matrix of Hermitian positive $(1,1)$ form $\sq\p\bp \phi$ on $ X$. By Lemma \ref{linear},
$$\left( \phi_{i\bar j}-\suml_{A,B=1}^{r-1} \phi^{A\bar B} \phi_{i\bar B} \phi_{A\bar
j}\right)$$ is a Hermitian positive $n\times n$ matrix. Since the
integrand is nonnegative, $R^{\tilde E}(u,u)=0$ if and only if \beq
\sum_{i,j}\sum_{\alpha,\beta } h_0( V_\alpha u_{}^{i\alpha})
\bar{(V_\beta u_{}^{j\beta})}\left( \phi_{i\bar
j}-\suml_{A,B=1}^{r-1}\phi^{A\bar B} \phi_{i\bar B} \phi_{A\bar
j}\right) \equiv 0\eeq on $X_s$ which means  $(u^{i\alpha})$ is a
zero matrix. In summary, we obtain
$$R^{\tilde E}(u,u)>0$$
for nonzero $u$, i.e. the induced metric $f$ on $\tilde E=S^k E\ts
F$ is Nakano-positive. \eproof \etheorem

 \bcorollary\label{Griffiths} If $E$ is ample, then for large $k$, $S^kE$
is Griffiths positive, i.e. there exists a Hermitian metric $h_k$ on
$S^kE$ such that $h_k$ is Griffiths-positive. \ecorollary

\subsection{Dual-Nakano-positivity}

By the curvature identity on $S^kE\ts F$,
$$ R_{i\bar j \alpha\bar \beta}=\int_{\ X_s} h_0
V_\alpha\bar V_\beta \phi_{i\bar j}\frac{\omega_s^{r-1}}{(r-1)!}
-\int_{X_s} h_0P_{i\alpha}\bar
P_{j\beta}\frac{\omega_s^{r-1}}{(r-1)!} $$ where $\phi$ is a weight
of the line bundle $\sO_{\P(E^*)}(k+r)\ts \pi^*(\det E^*)\ts
\pi^*(F)$. Although this line bundle can not be negative, it is
still possible that it is negative  in the local horizontal
direction, i.e. $(\phi_{i\bar j})$ is a Hermitian negative matrix.
For example, $F$ is a `` very negative" line bundle over $S$. If
$(\phi_{i\bar j})$ is Hermitian negative, then for any nonzero
$u=(u^{i\alpha})$, \be R_{i\bar j\alpha\bar \beta}
u_{}^{i\alpha}\bar {u_{}^{j\beta}}&=& \int_{X_s} h_0\phi_{i\bar j}(
V_\alpha u_{}^{i\alpha}) \bar{(V_\beta u_{}^{j\beta})}
\frac{\omega_s^{r-1}}{(r-1)!}\\&&-\int_{X_s} h_0P_{i\alpha}\bar
P_{j\beta}u^{i\alpha}\bar
u^{j\beta}\frac{\omega_s^{r-1}}{(r-1)!}\\&\leq& \int_{X_s}
h_0\phi_{i\bar j}( V_\alpha u_{}^{i\alpha}) \bar{(V_\beta
u_{}^{j\beta})} \frac{\omega_s^{r-1}}{(r-1)!}\\
&<&0 \ee Hence $S^kE\ts F$ is Nakano-negative. In the following, we
will prove that if $(S^{k+r}E\ts \det E^*\ts F)^*$ is ample, then
$S^kE\ts F$ is Nakano-negative which is equivalent to the statement:
if $S^{k+r}E\ts \det E^*\ts F$ is ample, then $S^kE\ts F$ is
dual-Nakano-positive. Here we use a well-known fact (\cite{D}):
\begin{quote}
\emph{ $E$ is dual-Nakano-positive if and only if $E^*$ is
Nakano-negative.}
\end{quote}

\noindent For simplicity, we  assume $k=1$ and $F=\det E$. In the
following we will show,  if $E^*$ is ample, then $E\ts \det E$ is
Nakano-negative.

 As similar as the quotient metric on
$\sO_{\P(E^*)}(1)$(see Proposition \ref{quotient} )  induced by the
morphism $(\pi^*E,\pi^*h)\>\sO_{\P(E^*)}(1)$, there is an induced
metric on $\sO_{\P(E)}(1)$ by the morphism
$(\pi^*(E^*),\pi^*h^*)\>\sO_{\P(E)}(1)$. For a fixed point $s\in S$,
we can choose a local  coordinate system $(z^1,\cdots, z^n)$
 and a local normal frame $(e_1,\cdots, e_r)$
of $E$ centered at point $s$. With respect to this trivialization,
we obtain:

\bproposition\label{Griffthsbase} If $(E,h)$ is Griffiths-positive,
then  the quotient metric $h^L$ on $L:=\sO_{\P(E)}(1)$ induced by
$(\pi^*E^*,\pi^*h^*)\>\sO_{\P(E)}(1)$  is negative  in the local
horizontal direction, i.e.  \beq \left(-\frac{\p^2\log h^L}{\p
z^i\p\bar z^j}\right) \eeq is Hermitian negative on the fiber
$X_s=\pi^{-1}(s)$ where $\pi:\P(E)\> S$.
 \bproof  %We can assume $(e_1,\cdots,
%e_r)$ is a normal frame of $E$ around $s\in S$.
Let $h_{\alpha\bar\beta}=h(e_\alpha,e_\beta)$ and $R_{i\bar
j\alpha\bar\beta}$ be the curvature components of $h$, then the
quotient metric on $\sO_{\P(E)}(1)$ is, \beq h^{L}=\frac{1}{\sum
h_{\alpha\bar\beta}W_\alpha\bar W_\beta}=\frac{1}{\sum
(\delta_{\alpha\beta}-R_{i\bar j \alpha\bar\beta}z^i\bar
z^j+O(|z|^3))W_\alpha\bar W_\beta} \eeq It is obvious that \beq
-\frac{\p^2\log h^L}{\p z^i\p\bar z^j}=-\sum_{\alpha,\beta}R_{i\bar
j \alpha\bar\beta}(s)\frac{W_\alpha\bar W_\beta}{|W|^2}\eeq which is
Hermitian negative on $X_s$ if $(E,h)$ is Griffiths-positive.
\eproof \eproposition

Let $v_k:E\>S^kE$ be the standard Veronese map which induces a map
\beq \bar v_k:\P(E)\>\P(S^kE) \eeq Let $\pi:\P(E)\>S$ and $\pi_k:
\P(S^kE)\>S$, then $\pi_k\circ \bar v_k=\pi$.  Now we fix a local
holomorphic coordinate system $(z^1,\cdots, z^n)$ centered at point
$s\in S$ and a local trivialization of $E$ and $S^kE$. It is obivous
that the map $\bar v_k$ sends $(z,W)$ to $(z,S^kW)$ where $S^kW$ is
the $k$-th symmetric power of homogeneous vector $W=[W_1,\cdots,
W_r]$, and so the horizontal part of $\bar v_k$ is identity. With
respect to this trivialization, we obtain

\btheorem\label{goodmetric} If $E$ is ample, then there exists a
Hermitian metric $h^{L}$ on $L=\sO_{\P(E)}(1)$ such that $h^L$ is
negative in  the horizontal direction, i.e. \beq
\left(-\frac{\p^2\log h^L}{\p z^i\p\bar z^j}\right) \eeq  is
Hermitian negative  on the fiber $X_s=\pi^{-1}(s)$ where
$\pi:\P(E)\> S$.\bproof By Corollary \ref{Griffiths}, for large $k$,
$S^kE$ is Griffiths-positive. By Proposition \ref{Griffthsbase},
there exists a Hermitian metric $\hat h_k$  on $\sO_{\P(S^kE)}(1)$,
such that $\hat h_k$ is Hermitian negative along the horizontal
direction. By the relation \beq \sO_{\P(E)}(k)=\bar
v_k^*\left(\sO_{\P(S^kE)}(1)\right) \eeq there is an induced metric
$h^L$ on $\sO_{\P(E)}(1)$  \beq h^L:=\left(\bar v_k^*(\hat
h_k)\right)^{\frac{1}{k}} \eeq Hence, we obtain \beq -\frac{\p^2\log
h^L}{\p z^i\p\bar z^j}=-\frac{1}{k}\frac{\p^2\log \hat h_k }{\p
z^i\p\bar z^j} \eeq since the horizontal direction of $\bar v_k$ is
 identity with respect to that trivialization.
 \eproof \etheorem

\btheorem\label{dualnakano} If $E^*$ is ample, then there exists a
Hermitian metric on $E\ts \det E$ which is  Nakano-negative.\bproof
By Theorem \ref{goodmetric}, if $E^*$ is ample, then there exists a
Hermitian metric $h^L$ on $L:=\sO_{\P(E^*)}(1)$ such that \beq
\left(-\frac{\p^2\log h^L}{\p z^i\p\bar z^j}\right) \eeq is
Hermitian negative. Let $\{\omega_s\}_{s\in S}$ be a smooth family
of Hermitian metric of the fiber $X_s$. We can set
$$h_0=\frac{(h^L)^{r+1}}{\det(\omega_s)}$$ and let
\beq \phi=-\log(h_0\det (\omega_s))=-(r+1)\log h^L \eeq Hence, we
obtain \beq \phi_{i\bar j}=-(r+1)\frac{\p^2\log h^L}{\p z^i\p\bar
z^j} \eeq Therefore $(\phi_{i\bar j})$ is Hermitian negative. On the
other hand, the metric induced by $h_0$ and $\{\omega_s\}_{s\in S}$
on $E\ts \det E$ has  curvature components \beq R_{i\bar j
\alpha\bar \beta}=\int_{ X_s} h_0 W_\alpha\bar W_\beta \phi_{i\bar
j}\frac{\omega_s^{r-1}}{(r-1)!} -\int_{X_s} h_0P_{i\alpha}\bar
P_{j\beta}\frac{\omega_s^{r-1}}{(r-1)!} \eeq Therefore, for any
nonzero $u=(u^{i\alpha})$, \be R_{i\bar j\alpha\bar \beta}
u_{}^{i\alpha}\bar {u_{}^{j\beta}} &\leq& \int_{X_s} h_0\phi_{i\bar
j}( W_\alpha u_{}^{i\alpha}) \bar{(W_\beta
u_{}^{j\beta})} \frac{\omega_s^{r-1}}{(r-1)!}\\
&<&0\ee The proof of  Nakano-negativity  of $E\ts \det E$ is
completed. \eproof \etheorem

Combined with Theorem \ref{Nakano}, Lemma \ref{positive} and Theorem
\ref{dualnakano} we obtain, \btheorem\label{main} Let $E$ be a
holomorphic vector bundle over a compact K\"ahler manifold $S$ and
$F$  a line bundle over $S$. Let $r$ be the rank of $E$ and $k\geq
0$ an arbitrary integer. If $S^{k+r}E\ts \det E^*\ts F$ is ample
over $S$, then $S^kE\ts F$ is both Nakano-positive and
dual-Nakano-positive.\etheorem

\subsection{Applications} \bcorollary\label{ske} If $E$ is an ample vector bundle
and $F$ is a nef line bundle, then there exists $k_0=k_0(S,E)$ such
that $S^{k}E\ts F$ is Nakano-positive and dual-Nakano-positive for
any $k\geq k_0$. In particular, $S^kE$ is Nakano-positive and
dual-Nakano-positive for $k\geq k_0$.  \bproof It is easy to see
that there exists  $k_0=k_0(S,E)$ such that for any $k\geq k_0$,
$S^{k+r}E\ts \det E^*$ is ample, and so is $S^{k+r}E\ts \det E^*\ts
F$. By Theorem \ref{main}, $S^kE\ts F$ is Nakano-positive and
dual-Nakano-positive. In particular, $S^kE$ is Nakano-positive and
dual-Nakano-positive for $k\geq k_0$.\eproof \ecorollary

\bcorollary\label{negative} If $E$ is an ample vector bundle and $F$
is a nef line bundle, or $E$ is a nef vector bundle and $F$ is an
ample line bundle,

\bd\item $S^{k}E\ts \det E\ts F$ is Nakano-positive and
dual-Nakano-positive for any $k\geq 0$.

\item If the rank  $r$ of $E$ is greater than $1$, then $S^mE^*\ts (\det
E)^t\ts F$ is Nakano-positive and dual-Nakano-positive if $t\geq
r+m-1$.

\ed

\bproof $(1)$ It follows by the ampleness of $S^{k+r}E\ts
F=S^{k+r}E\ts \det E^* \ts (\det E\ts F)$.

\noindent$(2)$ If $r>1$, it is easy to see $E^*\ts \det
E=\wedge^{r-1}E$. By the relation $$ S^{r+m}(E^*\ts\det E)\ts (\det
E)^{t-r-m+1} \ts F=S^{r+m}E^*\ts \det E\ts (\det E)^{t}\ts F$$ we
can apply Theorem \ref{main} to the pair $(E^*, (\det E)^t\ts F )$
and obtain the Nakano-positivity and dual-Nakano-positivity of
$S^mE^*\ts (\det E)^t\ts F$ when $t\geq r+m-1$. Let $E=T\P^2$, then
$E=E^*\ts \det E$ is Griffiths-positive but not Nakano-positive. So
we can not remove the restriction $t\geq r+m-1$. \eproof

\ecorollary

\bcorollary\label{optimal} If $S^{r+1}E\ts \det E^*$ is  ample, then
$E$ is Nakano-positive and dual-Nakano-positive and so $E$ is
Griffiths-positive. \ecorollary

\bremark  By Corollary \ref{optimal}, the ampleness of
$\sO_{\P(E^*)}(r+1)\ts \pi^*(\det E^*)$ implies the ampleness of
$\sO_{\P(E^*)}(1)$. But in general, the ampleness of
$\sO_{\P(E^*)}(1)$ can not imply the ampleness of
$\sO_{\P(E^*)}(r+1)\ts \pi^*(\det E^*)$.

\eremark

\section{Nakano-positivity and dual-Nakano-positivity of adjoint vector  bundles}

The following lemma is due to  Fujita (\cite{F1}) and Ye-Zhang
(\cite{YZ}). \blemma\label{fujita} Let $E$ be an ample vector bundle
over $S$. Let $r$ be the rank of $E$ and $n$ the dimension of $S$.
If $r\geq n+1$, then $\det E\ts K_S$ is ample except  $(S,E)\cong
(\P^n, \sO_{\P^n}(1)^{\ds {n+1}})$. \elemma

\btheorem\label{adjunction1} Let $E$ be an ample vector bundle over
$S$. Let $r$ be the rank of $E$ and $n$ the dimension of $S$.
\bd\item If $r>1$, then $S^{k}E\ts (\det E)^2\ts K_S$ is
Nakano-positive and dual-Nakano-positive for any $k\geq \max\{n-r,
0\}$.

\item If $r=1$, then the line bundle $E^{\ts (n+2)}\ts K_S$ is
Nakano-positive.
 \ed Moreover, the lower bound on $k$ is sharp. \bproof  $\bf{(1)}$ If $r>1$, then
$X=\P(E^*)$ is a $\P^{r-1}$ bundle which is not isomorphic to any
projective space. By Lemma \ref{fujita},
  $\sO_{\P(E^*)}(n+r)\ts K_X$ is ample. So
$$\sO_{\P(E^*)}(n)\ts \pi^*\left(K_S\ts\det E\right)$$
is ample and it is equivalent to the ampleness of $S^{n}E\ts (\det
E^*)\ts (\det E)^2\ts K_S$. If $k\geq \max\{n-r, 0\}$, $S^{r+k}E
\ts\det E^*\ts (\det E)^2\ts K_S$ is also ample, hence by Theorem
\ref{main}, $S^{k}E\ts (\det E)^2\ts K_S$ is Nakano-positive and
dual-Nakano-positive.

\noindent$\bf{(2)}$ It follows from Lemma \ref{fujita}. In fact, the
vector bundle $\tilde E=E^{\ds(n+2)}$ is an ample vector bundle of
rank $n+2$ and $\det\tilde E=E^{\ts (n+2)}$. By Lemma \ref{fujita},
 $\det\tilde E\ts K_S=E^{\ts(n+2)}\ts K_S$ is ample.

 Here
the lower bound $n-r$ is sharp.
 For any integer $k_0<n-r$, there exists some ample vector $E$ such that
 $E\ts (\det E)^{k_0}\ts K_S$ is not Nakano-positive, for example $(S,E)=(\P^4, \sO_{\P^4}(1)\ds
\sO_{\P^4}(1))$.
 \eproof \etheorem

\btheorem\label{adjunction2} Let $E$ be an ample vector bundle over
$S$. Let $r$ be the rank of $E$ and $n$ the dimension of $S$.  If
$r>1$, then $E\ts (\det E)^k\ts K_S$ is Nakano-positive and
dual-Nakano-positive for any $k\geq \max\{n+1-r, 2\}$. Moreover, the
lower bound is sharp. \bproof If $r\geq n-1$, by Theorem
\ref{adjunction1}, $E\ts \left(\det E\right)^2\ts K_S$ is
Nakano-positive and dual-Nakano-positive. Now we consider $1<r<n-1$.
By (\cite {I}, Theorem 2.5),  $K_S\ts \left(\det E\right)^{n-r}$ is
nef except the case $(S,E)=(\P^4, \sO_{\P^4}(1)\ds \sO_{\P^4}(1))$.
It is easy to check
$$S^{r+1}E\ts K_S\ts \left(\det
E\right)^{n-r}$$ is also ample in that case. By Theorem \ref{main},
$E\ts \left(\det E\right)^{n+1-r}\ts K_S$ is Nakano-positive and
dual-Nakano-positive. Here the lower bound $n+1-r$ is sharp. For any
integer $k_0<n+1-r$, there exists an ample vector bundle $E$ such
that $E\ts (\det E)^{k_0}\ts K_S$ is not Nakano-positive, for
example $(S,E)=(\P^4, \sO_{\P^4}(1)\ds \sO_{\P^4}(1))$. \eproof
\etheorem

\bremark In Theorem \ref{adjunction1} and \ref{adjunction2}, if
$r\geq n$, $E\ts (\det E)^2\ts K_S$ is Nakano-positive and
dual-Nakano-positive. If $E=T\P^n$, then $S^2E\ts \det E\ts
K_{\P^n}$ is Nakano-positive and dual-Nakano-positive.

\eremark

\noindent\bf{Problem:}  Is $S^2E\ts \det E\ts K_S$  Nakano-positive
and dual-Nakano-positive when $E$ is ample and $r\geq n$? If one can
show $S^{n+2}E\ts K_S$ is ample,
 or equivalently, $\sO_{\P(E^*)}(n+2)\ts\pi^*(K_S)$ is ample, by Theorem \ref{main}, $S^2E\ts \det E\ts K_S $
is Nakano-positive and dual-Nakano-positive.

\section{Vanishing theorems }

The following vanishing theorem is dual to Nakano(\cite{Na})(see
also Demailly(\cite{D})):
 \blemma\label{va} Let $E$ be a holomorphic vector bundle over a compact K\"ahler manifold $M$.
  If $E$ is Nakano-positive, then $H^{n,q}(M,E)=0$ for any
$q\geq 1$. If $E$ is dual-Nakano-positive, then $H^{q,n}(M,E)=0$ for
any $q\geq 1$.

\elemma

\btheorem\label{vanishing} Let $E, E_1,\cdots, E_\ell$ be
  vector bundles over an $n$-dimensional compact K\"ahler manifold $M$. Their ranks are  $r, r_1,\cdots,
r_\ell$ respectively. Let $L$ be a line bundle  on $M$. \bd

\item If $E$ is ample, $L$ is nef  and $r>1$, then
$$ H^{n,q}(M, S^kE\ts (\det E)^2\ts K_M\ts L)=H^{q,n}(M, S^kE\ts
(\det E)^2\ts K_M\ts L)=0 $$ for any $q\geq 1$ and $k\geq \max\{n-r,
0\}$.

\item If $E$ is ample, $L$ is nef  and $r>1$, then
$$ H^{n,q}(M, E\ts (\det E)^k\ts K_M \ts L)=H^{q,n}(M, E\ts (\det
E)^k\ts K_M \ts L)=0 $$ for any $q\geq 1$ and $k\geq \max \{n+1-r,
2\}$.

 \item  Let $r>1$. If $E$ is ample and $L$ is nef, or $E$ is
nef and $L$ is ample, then $$ H^{n,q}(M, S^mE^*\ts (\det E)^t\ts
L)=H^{q,n}(M, S^mE^*\ts (\det E)^t\ts L)=0 $$ for any $q\geq 1$ and
$t\geq r+m-1$.

\item If all $E_i$ are ample and $L$ is nef, or, all $E_i$ are nef and $L$ is ample,  then for any $k_1\geq 0,\cdots, k_\ell\geq 0$,
\be &&H^{n,q}(M,S^{k_1}E_1\ts\cdots \ts S^{k_\ell}E_\ell\ts \det
E_1\ts\cdots\ts \det E_\ell\ts L)\\&=&H^{q,n}(M,S^{k_1}E_1\ts\cdots
\ts S^{k_\ell}E_\ell\ts \det E_1\ts\cdots\ts \det E_\ell\ts L)=0 \ee
for $q\geq 1$.

\ed \bproof  By Theorem \ref{adjunction1}, Theorem \ref{adjunction2}
and Corollary \ref{negative},  the vector bundles in consideration
are all Nakano-positive and dual-Nakano-positive.  The results
follow from Lemma \ref{va}. \eproof
 \etheorem

\bremark  Part $(4)$ can be regarded as a generalization
 of Griffiths (\cite{G}, Theorem G) and
Demailly(\cite{D1}, Theorem 0.2).

 \eremark

The following results generalize  Griffiths' vanishing theorem( see
also \cite{LN2}, Corollary 1.5):

\bproposition\label{induction} Let $r$ be the rank of $E$ and $k\geq
1$. For any $t\geq 0$, if $S^{t+kr}E\ts L$ is ample,
$$H^{n,q}(M,S^{t}E\ts (\det E)^{k}\ts L)=H^{q,n}(M,S^{t}E\ts (\det E)^{k}\ts L)=0$$ for any $q\geq 1$.
\bproof By Theorem \ref{main},  $S^tE\ts (\det E)^k\ts L$ is
Nakano-positive and dual-Nakano-positive. The results follow by
Nakano's vanishing theorem. \eproof \eproposition \bremark Theorem
\ref{main0} allows us to do induction to deduce more positivity
results. For example, if $S^{m}E\ts L$ is ample, then $S^{m-r}E\ts
\det E\ts L$ is (dual-)Nakano-positive and so it is ample. Using
Theorem \ref{main0} again, we get $S^{m-2r}\ts (\det E)^2\ts L$ is
Nakano-positive and dual-Nakano-positive. Finally, we get $S^{t}E\ts
(\det E)^k\ts L$ is Nakano-positive and dual-Nakano-positive, if
$m=t+kr$ for some $0\leq t<r$. It is obvious that the
(dual-)Nakano-positivity turns stronger and stronger under
induction. This explains why a lot of vanishing theorems involve a
power of $\det E$. \eremark

If $L$ is an ample line bundle over a compact K\"ahler manifold $M$
and $F$ is an arbitrary line bundle over $M$. By comparing the Chern
classes,  there exists a constant $m_0$ such that $L^{m_0}\ts F$ is
ample and so it is positive. If $E$ is an ample vector bundle and
$F$ is an arbitrary vector bundle, it is easy to see $S^kE\ts F$ is
ample for large $k$. But, in general, we don't know whether an ample
vector bundle carries a Griffiths-positive or Nakano-positive
metric. In the following, we will construct Nakano-positive and
dual-Nakano-positive metrics on various  ample vector bundles.

\blemma\label{power1} If $L$ is an ample line bundle over $M$ and
$F$ is an arbitrary vector bundle. There exists an integer $m_0$
such that $L^{m_0}\ts F$ is Nakano-positive and
dual-Nakano-positive. \bproof Let $h_0$ be a positive metric on $L$
and $\omega$ be the curvature of $h_0$ which is also the K\"ahler
metric fixed on $M$. For any metric $g$ on $F$, the curvature $R^g$
has a lower bound in the sense \beq \min_{x\in M}\inf_{u\neq 0}
\frac{R^g(u(x),u(x))}{|u(x)|^2}\geq -(m_0-1) \eeq where $u\in
\Gamma(M,T^{1,0}M\ts F)$. The curvature of  metric $h^{m_0}\ts g$ on
$L^{m_0}\ts F$ is given by \beq \hat R=m_0\omega\cdot g+h_0^m\cdot
R^g \eeq Therefore
$$\hat R(v\ts u,v\ts u)\geq |u|^2 h_0^{m_0}(v,v)$$
for any $v\in \Gamma(M,L^{m_0})$ and $u\in\Gamma(M,T^{1,0}M\ts F)$.
\eproof \elemma

\blemma\label{power2} If $E$ is (dual-)Nakano-positive and $F$ is a
nef \emph{line} bundle, then $E\ts F$ is (dual-)Nakano-positive.
\bproof Fix a K\"ahler metric on $M$. Let $g$ be a Nakano-positive
metric on $E$, then there exists $2\eps>0$ such that
$$R^g(u(x),u(x))\geq 2\eps |u(x)|^2$$
for any $u\in\Gamma(M,T^{1,0}M\ts E)$. On the other hand, by a
result of \cite{DPS}, there exists a smooth metric $h_0$ on the nef
line bundle $F$ such that \beq R^{h_0}\geq -\eps \omega h_0 \eeq The
curvature of $g\ts h_0$ on $E\ts F$ is
$$\hat R= R^g\cdot h_0+ g\cdot R^{h_0}$$
For any $u\in \Gamma(M,T^{1,0}M\ts E)$ and $v\in \Gamma(M,F)$ \beq
\hat R(u\ts v,u\ts v)\geq \left(R^{g}(u,u)-\eps |u|^2\right)
h_0(v,v)\geq \eps |u|^2h_0(v,v) \eeq For  dual-Nakano-positivity,
the proof is similar.
 \eproof \elemma

\btheorem\label{nadu} If $E$ is an ample vector bundle and $F$ is an
arbitrary vector bundle over $M$, then there exists $k_0=k_0(M,E,F)$
such that $S^kE\ts F$ is Nakano-positive and dual-Nakano-positive
for any $k\geq k_0$. \bproof By Lemma \ref{power1}, there exists
$m_0$ such that $(\det E)^{m_0}\ts F$ is Nakano-positive and
dual-Nakano-positive. On the other hand, there exists
$k_0=k_{0}(E,m_0,M)$ such that $\sO_{\P(E^*)}(r+k)\ts \pi^*(\det
E^*)^{m_0+1}$ is ample for $k\geq k_0$. It is equivalent to the
ampleness of vector bundle $S^{r+k}E\ts (\det E^*)^{m_0+1}$. By
Theorem \ref{main}, $S^{k}E\ts (\det E^*)^{m_0}$ is Nakano-positive
and dual-Nakano-positive. Since the tensor product of two
(dual-)Nakano-positive vector bundles  is (dual-)Nakano-positive,
$S^kE\ts F= (S^kE \ts (\det E^*)^{m_0}) \ts \left((\det E)^{m_0}\ts
F\right)$ is Nakano-positive and dual-Nakano-positive for $k\geq
k_0$. \eproof \etheorem

\noindent The following results are well-known in algebraic
geometry, but merit a proof in our setting.

 \bcorollary\label{arbitrary} If $E$ is ample over $M$, $L$
is a nef line bundle and $F$ is an arbitrary vector bundle, \bd\item
there exists $k_0=k_0(M,E,F)$ such that for any $k\geq k_0$.$$
H^{p,q}(M,S^kE\ts F)=0$$ for $q\geq 1$ and $p\geq 0$.
\item
 there exists
$k_0=k_0(M,E)$ such that for any $k\geq k_0$, $$ H^{p,q}(M,S^kE \ts
L)=0 $$ for any $q\geq 1$ and $p\geq 0$.\ed \bproof $\bf{(1)}$ By
Theorem \ref{nadu}, there exists $k_0=k_0(M,E,F)$ such that $S^kE\ts
F\ts \Lambda^{n-p}T^{1,0}M$ is Nakano-positive for any $p$. On the
other hand $$ H^{p,q}(M,S^kE\ts F)=H^{n,q}(M,S^kE\ts F\ts
\Lambda^{n-p}T^{1,0}M) $$  By Nakano vanishing theorem,   $
H^{p,q}(M,S^kE\ts F)=0$ for $q\geq 1$ and $p\geq 0$ if $k\geq k_0$.
The proof of part $\bf{(2)}$ is similar. \eproof \ecorollary

\section{Comparison of Griffiths-positive and Nakano-positive metrics}

Let $(E, h)$ be a Hermitian vector bundle. In general, it is not so
easy to write down the exact curvature formula of $(S^kE, S^kh)$. In
this section, we give an algorithm to compute the curvature of
$(S^kE, S^kh)$. As applications, we can disprove  the
Griffiths-positivity and Nakano-positivity of a given metric on
$\P^n$.

 Let $h$ be a Hermitian metric
on $E$, $h^L$ be the induced metric  in (\ref{inducedmetric}) on
$L=\sO_{\P(E^*)}(1)$. Let $F$ be a line bundle with Hermitian metric
$h^F$. Naturally, there is an induced metric $S^kh\ts h^F$ on the
vector bundle $S^kE\ts F$. On the other hand, we can construct a new
metric $f$ on $S^kE\ts F$ by  formula (\ref{l2metric}).  There is a
\bf{canonical way} to do it. Let ${\tilde L}=L^k \otimes \pi^*(F)$.
The induced metric on $\tilde L$
  is $h_0=(h^L)^k\ts \pi^*(h^F)$ and the induced metric on   $\det(T_{X/S})=L^r\ts \pi^*(\det E^*)$  is $(h^L)^r\ts
\pi^*(  \det(h)^{-1})$. These two metrics induce a metric
$\lambda_0=(h^L)^{k+r}\ts
 \pi^*\left( h^F\cdot \det(h)^{-1}\right)$ on $\tilde L\ts \det (T_{X/S})$. Now we can
  polarize each fiber $X_s$ by the curvature of $\lambda_0$.
 By formula (\ref{inducedcurvature}),
\beq \omega_s=-\frac{\sq }{2\pi}\p_s\bp_s\log
\lambda_0=\frac{(k+r)\sq }{2\pi}\p_s\bp_s\log\left(\sum
h^{\alpha\bar\beta}W_\alpha\bar
W_\beta\right)=(k+r)\omega_{FS}\label{43}\eeq By a simple linear
algebraic argument, we obtain\beq
\frac{\lambda_0}{\det(\omega_s)}=\frac{(h^L)^k\ts \pi^*(
h^F)}{(k+r)^{r-1}}=\frac{h_0}{(k+r)^{r-1}} \label{44}\eeq
 Now we can use $(\tilde L, h_0)$ and $(X_s, \omega_s)$ to construct a ``new" metric $f$ on $S^kE\ts F$ by formula (\ref{l2metric}).

\btheorem\label{invariant} The metric $f$  has the form\beq
f=\frac{(r+k)^{r-1}}{(r+k-1)!}\cdot S^k h\ts h^F\eeq Moreover,  $f$
is a constant multiple of the metric constructed in Theorem
\ref{main}.\bproof Without loss of generality, we can choose  normal
coordinates for the metric $h$ at a fix point $s\in S$. By formula
(\ref{inducedcurvature}), the metric $h_0=(h^L)^k\ts h^F$ on $L^k\ts
F$ induced by $(E,h)$ and $(F,h^F)$ can be written as
$\frac{h^F}{|W|^{2k}}$ locally on the fiber $X_s\cong\P^{r-1}$. By
formula (\ref{43}), the metric $f$ defined by (\ref{l2metric}) has
the following form
$$f_{\alpha\bar \beta}=\int_{X_s} h_0 V_\alpha\bar V_\beta \frac{\omega_s^{r-1}}{(r-1)!}=(k+r)^{r-1} h^F\int_{\P^{r-1}} \frac{V_\alpha\bar V_\beta}{|W|^{2k}}\frac{\omega_{FS}^{r-1}}{(r-1)!}$$
Here $V_\alpha, V_\beta$ are homogeneous monomials of degree $k$ in
$W_1,\cdots ,W_r$. By Lemma \ref{calculation},
 $$f_{\alpha\beta}=\frac{(r+k)^{r-1}}{(r+k-1)!} \delta_{\alpha\beta} h^F$$
 that is $f=\frac{(r+k)^{r-1}}{(r+k-1)!}\cdot S^k
h\ts h^F$. By  formulas (\ref{44}) and (\ref{45}),  $f$ is a
constant multiple of the metric constructed in Theorem \ref{main}.
 \eproof \etheorem

 \btheorem\label{DS} If $(E,h)$ is a Griffiths-positive vector bundle, then
 \bd\item $(S^kE \ts (\det E)^\ell, S^kh \ts (\det h)^\ell)$ is
 Nakano-positive and dual-Nakano-positive for any $k\geq 0$ and $\ell\geq
 1$.
\item There exists $k_0=k_0(M,E)$ such that $(S^kE,S^kh)$ is
Nakano-positive and dual-Nakano-positive for any $k\geq k_0$. \ed

\bproof These  follow by Theorem \ref{main} and Theorem
\ref{invariant}. \eproof
 \etheorem

\bcorollary\label{seminakano} Let $h_{FS}$ be the Fubini-Study
metric on $T\P^n$ with $n\geq 2$, then
 \bd\item $(S^{n+1}T\P^n\ts K_{\P^n}, S^{n+1}h_{FS}\ts \det(h_{FS})^{-1})$ is
semi-Griffiths-positive. Moreover, $S^{n+1}T\P^n\ts K_{\P^n}$ can
not admit a Griffiths-positive metric.

 \item $(T\P^n, h_{FS})$
is dual-Nakano-positive and semi-Nakano-positive.

\item $(S^{k}T\P^n\ts K_{\P^n}, S^{k}h_{FS}\ts
\det(h_{FS})^{-1})$ is Griffiths-positive for any $k\geq n+2$.

\item $(S^{k}T\P^n,S^kh_{FS})$ is Nakano-positive and dual-Nakano-positive for any $k\geq 2$.
 \ed
\bproof \bf{(1)} By the Euler  sequence \beq 0\>\C\>T\P^n\>
\sO_{\P^n}(1)^{\ds (n+1)}\>0 \label{el}\eeq we know $T\P^n\ts
\sO_{\P^n}(-1)$ is the quotient bundle of trivial bundle
$\C^{\ds(n+1)}$. Hence
$$S^{n+1}T\P^n\ts
K_{\P^n}=S^{n+1}\left(T\P^n\ts\sO_{\P^n}(-1)\right)$$ with the
canonical metric is semi-Griffiths-positive.
 However, if
$S^{n+1}T\P^n\ts K_{\P^n}$ admits a Griffiths-positive metric, by
Corollary \ref{optimal}, $T\P^n$ is Nakano-positive which is
impossible for $n\geq 2$. \noindent$\bf{(2)}$  The curvature of
$E=T\P^n$ with respect to the standard Fubini-Study metric $h_{FS}$
is \beq R_{i\bar j k\bar \ell}=h_{i\bar j}h_{k\bar
\ell}+h_{i\bar\ell}h_{k\bar j} \eeq Without loss of generality, we
assume $h_{i\bar j}=\delta_{ij}$ at a fixed point, then  \beq
R_{i\bar j k\bar \ell} u^{ik}\bar u^{j\ell}=\frac{1}{2}\sum_{j,k}
|u^{jk}+u^{kj}|^2\label{pn}\eeq which means that $(E,h_{FS})$ is
semi-Nakano-positive but not Nakano-positive. For the
dual-Nakano-positivity of $(T\P^n,h_{FS})$ we can check that by
definition. We can also show it by the monotone property of
dual-Nakano-positivity of quotient bundles. By the Euler sequence
(\ref{el}), $T\P^n$ is the quotient bundle of dual-Nakano-positive
bundle $\sO_{\P^n}(1)^{\ds (n+1)}$ and so $T\P^n$ is
dual-Nakano-positive.

\noindent $\bf{(3)}$ It follows by the identity
$$S^kT\P^n\ts K_{\P^n}=S^{k}(T\P^n\ts \sO_{\P^{n}}(-1))\ts
\sO_{\P^n}(k-n-1)$$ and semi-Griffiths positivity of $T\P^n\ts
\sO_{\P^n}(-1)$.

\noindent $\bf{(4)}$ By Theorem \ref{main},  the canonically induced
metric $f$ is Nakano-positive and dual-Nakano-positive. On the other
hand, by Theorem \ref{invariant},  $f$ is a constant multiply of
$S^kh_{FS}$. The lower bound of $k$ follows from $\bf{(1)}$ and
$\bf{(2)}$. \eproof \ecorollary

\bproposition\label{good} \bd\item $(E,h)$ is Griffiths-positive if
and only if $(S^kE,S^kh)$ is Griffiths-positive for some $k\geq 1$.

\item If $(E,h)$ is (dual-)Nakano-positive, then $(S^kE,S^kh)$
is (dual-)Nakano-positive for any $k\geq 1$.\ed \bproof By Theorem
\ref{invariant}, $S^kh$ is a constant multiple of the metric
constructed by formula (\ref{l2metric}). So by Theorem
\ref{curvature}, we can write down the curvature formula of $S^kh$
explicitly. In a normal coordinates of $h$ at a fixed point, the
curvature formula (\ref{curvatureformula}) can be simplified by
Lemma \ref{calculation}. We obtain curvature formulas (\ref{curv2})
and (\ref{curv8}).

For the convenience of the reader, we assume $k=2$ at first. We can
choose
 normal coordinates at a fixed point. Let $\{e_1,\cdots, e_r\}$
be the local basis at that point. The ordered basis of $S^2E$ at
that point are $\{e_{1}\ts e_1, e_{1}\ts e_2,\cdots,e_r\ts e_{r-1},
e_{r}\ts e_r\}$. We  denote them by $e_{(\alpha,\beta)}=e_\alpha\ts
e_\beta$ with $\alpha\leq \beta$. The curvature tensor $S^2h$ is
\beq R_{i\bar j(\alpha,\gamma)\bar{(\beta,\delta)}}=R_{i\bar j
\alpha\bar\beta} \delta_{\gamma\delta}+R_{i\bar j
\gamma\bar\delta}\delta_{\alpha\beta}+R_{i\bar j \gamma\bar
\beta}\delta_{\alpha\delta}+R_{i\bar j \alpha\bar
\delta}\delta_{\gamma\beta}\label{curv2}\eeq where $R_{i\bar
j\alpha\bar \beta}$ is the curvature tensor of $E$.  Let
$u=\suml_i\suml_{\alpha\leq \gamma}u_{i(\alpha,\gamma)}
e_{(\alpha,\gamma)}\in \Gamma(M,T^{1,0}M\ts S^2E)$. For simplicity
of notations, we  extend the values of $u_{i(\alpha,\gamma)}$ to all
indices $(\alpha,\gamma)$ by  setting $u_{i(\alpha,\gamma)}=0$ if
$\gamma<\alpha$. Therefore
\begin{eqnarray}&&\nonumber\sum_{i,j}\sum_{\stackrel{\alpha\leq \gamma}{ \beta\leq
\delta}} R_{i\bar j(\alpha,\gamma)\bar{(\beta,\delta)}}
u_{i(\alpha,\gamma)}\bar{u}_{j(\beta,\delta)}\\&=&\nonumber
\sum_{i,j}\sum_{{\alpha, \gamma,}{ \beta, \delta}} R_{i\bar
j(\alpha,\gamma)\bar{(\beta,\delta)}}
u_{i(\alpha,\gamma)}\bar{u}_{j(\beta,\delta)}\\
&=&\nonumber\sum_{i,j,\alpha,\beta,\gamma,\delta} \big( R_{i\bar
j\alpha\bar\beta} u_{i(\alpha,\gamma)}\bar
u_{j(\beta,\gamma)}+R_{i\bar j
\gamma\bar\delta}u_{i(\alpha,\gamma)}\bar
u_{j(\alpha,\delta)}\\\label{skcurvature}&&+R_{i\bar j
\gamma\bar\beta}u_{i(\alpha,\gamma)}\bar
u_{j(\beta,\alpha)}+R_{i\bar j \alpha\bar\delta}
u_{i(\alpha,\gamma)}\bar u_{j(\gamma,\delta)}\big)\\\nonumber &=&
\sum_{\gamma}\sum_{i,j,\alpha,\beta} R_{i\bar j\alpha \bar
\beta}\left(u_{i(\alpha,\gamma)}+u_{i(\gamma,\alpha)}\right)\bar{\left(
u_{j(\beta,\gamma)}+u_{j(\gamma,\beta)}\right)}\end{eqnarray} Hence
$(S^2E, S^2 h)$ is Nakano-positive if $(E,h)$ is Nakano-positive.
For the general case, we set $A=(\alpha_1,\cdots, \alpha_k)$ and
$B=(\beta_1,\cdots, \beta_k)$ with $\alpha_1\leq \cdots \leq
\alpha_k$ and $\beta_1\leq \cdots\leq \beta_k$. The basis of $S^kE$
are $\{e_{A}=e_{\alpha_1}\ts \cdots \ts e_{\alpha_k}\}$.  The
curvature tensor of $(S^kE, S^kh)$ is \beq R_{i\bar j A\bar
B}=\sum_{\alpha,\beta=1}^r\sum_{s,t=1}^kR_{i\bar
j\alpha\bar\beta}\delta_{\alpha
\alpha_s}\delta_{\beta\beta_t}\delta_{A_sB_t} \label{curv8} \eeq
where $A_s=(\alpha_1,\cdots, \alpha_{s-1}, \alpha_{s+1},\cdots,
\alpha_k)$, $B_t=(\beta_1,\cdots,\beta_{t-1},\beta_{t+1},\cdots,
\beta_{k})$
 and $\delta_{A_sB_t}$ is the multi-index delta function( see formula (\ref{generalizeddelta})).  We have the curvature formula,
\begin{eqnarray}&&\nonumber\sum_{i,j, A, B} R_{i\bar j A\bar B} u_{iA}\bar u_{j
B}\\&=&\sum_{\alpha_{1},\cdots,\alpha_{k-1}} \sum_{\sigma\in
S_{k-1}}\sum_{i,j,\alpha,\beta} R_{i\bar j\alpha\bar \beta}
V_{i\alpha\alpha_{\sigma(1)}\cdots \alpha_{\sigma(k-1)}}\bar
V_{j\beta\alpha_{\sigma(1)}\cdots \alpha_{\sigma(k-1)}}\label{dd}
\end{eqnarray} where $S_{k-1}$ is the permutation group in $(k-1)$ symbols and
$$V_{i\alpha \alpha_1\cdots \alpha_{k-1}}=\sum_{s=1}^k u_{i A^s},\ \
\ A^s=(\alpha_1,\cdots, \alpha_{s-1},\alpha,\alpha_{s+1},\cdots,
\alpha_k)$$  The Nakano-positivity of $(S^kE,S^kh)$ follows
immediately from the Nakano-positivity of $(E,h)$ by  formula
(\ref{dd}). With the help of curvature formula (\ref{curv8}), we can
prove Griffiths-positivity and dual-Nakano-positivity  of $S^kE$ in
a similar way. Here, we  use another way to show it. $S^kE$ can be
viewed as a quotient bundle of $E^{\ts k}$. If $(E,h)$ is
Griffiths-positive(resp. dual-Nakano-positive), $(E^{\ts k},h^{\ts
k})$ is Griffiths-positive(resp. dual-Nakano-positive) and so the
quotient bundle $S^kE$ is Griffiths-positive(resp.
dual-Nakano-positive)(\cite{D}). The induced metrics on quotient
bundles are exactly the given ones.
 \eproof \eproposition

\bremark Part $(1)$ is an analogue of ampleness: $E$ is ample if and
only if $S^kE$ is ample for some $k\geq 1$.
 The converse of part $(2)$ is not valid in general. We know $(S^2
T\P^n,S^2 h_{FS})$ is Nakano-positive, but $(T\P^n, h_{FS})$ is not
Nakano-positive as shown in the following. \eremark

\bexample In this example, we will show the Nakano-positivity of
$(S^2T\P^2, S^2 h_{FS})$ in local coordinates. At a fixed point, we
choose a normal coordinates of $T\P^2$. Let $\{e_1,e_2\}$ be the
ordered basis of $T\P^2$ at that point. The ordered basis of
$S^2T\P^2$ are $e_{(1,1)}=e_1\ts e_1, e_{(1,2)}=e_1\ts e_2$ and
$e_{(2,2)}=e_2\ts e_2$. Using the same notation as Proposition
\ref{good}, we set
$V_{i\alpha\gamma}=u_{i(\alpha,\gamma)}+u_{i(\gamma,\alpha)}$ where
$u=\suml_i\suml_{\alpha\leq \gamma}u_{i(\alpha,\gamma)} \frac{\p}{\p
z^i}\ts e_{(\alpha,\gamma)}\in \Gamma(\P^2, T^{1,0}\P^2\ts
S^2T\P^2)$. For $\gamma=1$, the $2\times 2$ matrix
$\left(V_{i\alpha1}\right)$ has the form
$$
T_1=\left(%
\begin{array}{cc}
 2u_{1(1,1)} &u_{1(1,2)} \\
  2u_{2(1,1)} & u_{2(1,2)} \\
\end{array}%
\right) $$ For $\gamma=2$,  the $2\times 2$ matrix
$\left(V_{i\alpha2}\right)$ is
$$
T_2=\left(%
\begin{array}{cc}
 u_{1(1,2)} &2u_{1(2,2)} \\
  u_{2(1,2)} & 2u_{2(2,2)} \\
\end{array}%
\right) $$ The total $2\times 3$ matrix
$\left(u_{i(\alpha,\beta)}\right)$ is
$$
T=\left(%
\begin{array}{ccc}
 u_{1(1,1)} &u_{1(1,2)}& u_{1(2,2)} \\
  u_{2(1,1)} & u_{2(1,2)}& u_{2(2,2)} \\
\end{array}%
\right) $$ By formulas (\ref{skcurvature}) and (\ref{pn}), \be
\sum_{i,j,\alpha, \gamma, \beta, \delta} R_{i\bar
j(\alpha,\gamma)\bar{(\beta,\delta)}}
u_{i(\alpha,\gamma)}\bar{u}_{j(\beta,\delta)}&=&\sum_{i,j,\alpha,
\beta}\left(R_{i\bar j\alpha\bar\beta}V_{i\alpha1}\bar V_{j\beta 1}+
R_{i\bar j\alpha\bar\beta}V_{i\alpha2}\bar V_{j\beta
2}\right)\\&=&\frac{1}{2}\sum_{i,\alpha}|V_{i\alpha 1}+V_{\alpha i
1}|^2+\frac{1}{2}\sum_{i,\alpha}|V_{i\alpha 2}+V_{\alpha i2}|^2\ee
It equals zero if and only if $T_1$ and $T_2$ are skew-symmetric
which means $T\equiv 0$. The Nakano-positivity of  $(S^2T\P^2,
S^2h_{FS})$ is proved.

\eexample

\bibliographystyle{amsplain}

\begin{thebibliography}{10}

\bibitem{BS} M.C. Beltrametti and  A.J. Sommese, \textit{The adjunction theory of
complex projective varieties}.  Expositions in Math. 16, de Gruyter,
1995.

\bibitem{Bo1} B. Berndtsson, \textit{Curvature of vector bundles associated to holomorphic fibrations}. Ann. of Math. (2) \bf{169} (2009), no. 2, 531--560

\bibitem{Bo2} B. Berndtsson,  \textit{Positivity of direct image bundles and convexity on the space of K\"ahler metrics}. J. Differential Geom. \bf{81}(2009), no. 3, 457--482.

\bibitem{Bo3} B. Berndtsson, \textit{Strict and non strict positivity of direct image bundles}.  arXiv:1002.4797.

\bibitem{Be1}  R. Berman,  \textit{Bergman kernels and equilibrium measures for line bundles over projective manifolds}. Amer. J. Math. \bf{131} (2009), no. 5, 1485--1524.

\bibitem{Be2} R. Berman, \textit{Relative K\"ahler-Ricci flows and their quantization}.  arXiv:1002.3717.

\bibitem{CF}  F. Campana and H. Flenner, \textit{A characterization of ample vector bundles on a curve}. Math. Ann. \bf{287} (1990), no. 4, 571--575.

\bibitem{CH}  M. Cornalba and J. Harris,  \textit{Divisor classes associated to families of stable varieties, with applications to the moduli space of curves}. Ann. Sci. \'ecole Norm. Sup. (4) \bf{21} (1988), no. 3, 455--475.

\bibitem{D} J-P. Demailly,
\textit{Complex analytic and algebraic geometry}. book online
http://www-fourier.ujf-grenoble.fr/~demailly/books.html.

\bibitem{D1} J-P. Demailly, \textit{Vanishing theorems for tensor powers of an ample vector bundle}. Invent.
Math. \bf{91} (1988), 203--220.



\bibitem{DPS} J-P. Demailly, T. Peternell and M. Sehneider,  \textit{Compact complex manifolds with numerically
effective tangent bundles}. J. Alg. Geom. \bf{3}(1994), 295--345

\bibitem{DS} J-P. Demailly and H. Skoda,  \textit{Relations entre les notions de positivit\'{e} de P.A. Griffiths et de
S. Nakano}, S\'eminaire P. Lelong-H. Skoda (Analyse), ann\'ee
1978/79, Lecture notes in Math. No. 822, Springer-Verlag, Berlin
(1980) 304--309.



\bibitem{F1} T. Fujita, \textit{ On adjoint bundles of ample vector bundles}. Proc. Complex Alg.
Var., Bayreuth 1990 ; Lecture Notes in Math. 105--112. Springer
1991.


\bibitem{G} P. Griffiths, \textit{ Hermitian differential geometry, Chern classes and positive vector bundles},
Global Analysis, papers in honor of K. Kodaira, Princeton Univ.
Press, Princeton (1969) 181--251.

\bibitem{GH} P. Griffiths and J. Harris,  \textit{Principles of algebraic geometry}. Wiley, New York (1978).

\bibitem{Har} R. Hartshorne,\textit{Ample vector bundles}, Publ. Math. I.H.E.S. \bf{29} (1966) 319--350.

\bibitem{H2} L. H\"ormander,  \textit{$L^2$ estimates and existence theorems for the $\bp$ operator}, Acta Math. \bf{113}
(1965) 89--152.

\bibitem{H} L. H\"ormander, \textit{ An introduction to complex analysis in several complex variables}. D. Van Nostrand Co., Inc., 1966.

\bibitem{I} H. Ishihara, \textit{Some adjunction properties of ample vector
bundles}. Canad. Math. Bull. \bf{44} (2001) No 4. 452--458.


\bibitem{L} R. Lazarsfeld,
\textit{Positivity in algebraic geometry.}  I, II.  Ergebnisse der
Mathematik und ihrer Grenzgebiete. 3. Folge. A Series of Modern
Surveys in Mathematics  Springer-Verlag, Berlin, 2004.



\bibitem{LN} F. Laytimi and  W. Nahm, \textit{ A generalization of Le Potier's vanishing theorem},
Manuscripta math. \bf{113} (2004), 165--189.

\bibitem{LN2} F. Laytimi and  W. Nahm, \textit{ On a vanishing problem of Demailly}. Int.
Math. Res. Not.  \bf{47}(2005), 2877--2889.

\bibitem{LSY} K. Liu,  X. Sun and  S-T. Yau, \textit{Canonical metrics on the moduli space of Riemann Surfaces
I}.  J. Differential Geom. \bf{68} (2004), 571--637.

\bibitem{MZ} X. Ma and W. Zhang, \textit{Superconnection and family Bergman kernels}. C. R. Math. Acad. Sci. Paris \bf{344} (2007), no. 1,
41--44.

\bibitem{M} L. Manivel, \textit{ Vanishing theorems for ample vector bundles}. Invent. math. \bf{127}(1997),
401--416.


\bibitem{Mo}  S. Mori, \textit{ Projective manifolds with ample tangent bundles}. Ann. of Math. (2) \bf{110} (1979), no. 3, 593--606.

\bibitem{MT1} C. Mourougane and  S. Takayama, \textit{ Hodge metrics and positivity of direct images}. J. Reine Angew. Math. \bf{606} (2007), 167--178.

\bibitem{MT2} C. Mourougane and  S. Takayama, \textit{ Hodge metrics and the curvature of higher direct images}. Ann. Sci. \'Ec. Norm. Sup\'er. (4) \bf{41} (2008) no. 6, 905--924.


\bibitem{Na} S. Nakano,\textit{  On complex analytic vector bundles}, J. Math. Soc. Japan \bf{7} (1955) 1--12.


\bibitem{P} T. Peternell, \textit{ A characterisation of $\P^n$ by vector bundles}. Math. Z.
\bf{205}(1990), 487--490.

\bibitem{PLS} T. Peternell, J. Le Potier and M. Schneider, \textit{Vanishing theorems,
linear and quadratic normality}. Invent. Math. \bf{87} (1987),
573--586.
\bibitem{S1} G. Schumacher, \textit{ On the geometry of moduli spaces}. Manuscripta
Math. \bf{50} (1985), 229--267.


\bibitem{S} G. Schumacher, \textit{  Curvature of higher direct images and applications}. arXiv:1002.4858.


\bibitem{SS} B. Shiffman and A.J. Sommese,  \textit{Vanishing theorems on complex manifolds}. Progress in Mathematics, \bf{56}. Birkh\"auser 1985.

\bibitem{Siu} Y.T. Siu, \textit{Curvature of the Weil-Petersson metric in the moduli space of compact K\"ahler-Einstein manifolds of negative first
Chern class}. {Contributions to Complex Analysis}, Papers in Honour
of Wilhelm Stoll. Vieweg, Braunschweig, 1986.

\bibitem{SY} Y.T. Siu and S-T. Yau, \textit{ Compact K\"ahler manifolds of positive bisectional curvature}. Invent. Math. \bf{59} (1980), no. 2, 189--204.


\bibitem{T} G. Tian,  \textit{On a set of polarized K\"ahler metrics on algebraic manifolds}. J. Differential Geom. \bf{32} (1990), no. 1, 99--130.

\bibitem{U}  H. Umemura, \textit{ Some results in the theory of vector bundles}. Nagoya Math. J. \bf{52} (1973), 97--128.



\bibitem{Wo}  S. Wolpert, \textit{ Chern forms and the Riemann tensor for the moduli space of curves}. Invent. Math. \bf{85} (1986), no. 1, 119--145.


\bibitem{YZ} Y-G. Ye and Q. Zhang, \textit{ On ample vector bundles whose adjunction bundles are not numerically
effective}. Duke Math. J. \bf{60} (1990), 671--687.







{
\bigskip\noindent  $^{\dagger}$ \textsf{ Department of Mathematics,  University of California at Los Angeles,\\
\ \ \ Los Angeles, CA,  90095-1555}\\
 \emph{E-mail Address}: liu@math.ucla.edu;\ \ \ xkyang@math.ucla.edu


\bigskip \noindent $^*$ \textsf{Department of Mathematics,  Lehigh University, Bethlehem, PA
18015}\\ \emph{E-mail Address}: xis205@lehigh.edu}
\\
\\


\end{thebibliography}

\end{document}